\pdfoutput=1

\documentclass[a4paper, 11pt, abstracton, titlepage, oneside]{scrartcl}
\makeatletter

\usepackage[T1]{fontenc}
\usepackage[utf8]{inputenc}
\usepackage[onehalfspacing]{setspace}
\usepackage[headsepline]{scrlayer-scrpage}
\clearpairofpagestyles 
\usepackage{layout}
\usepackage{geometry}
\usepackage{adjustbox}
\usepackage{authblk}
\usepackage[titletoc,title]{appendix}

\usepackage[hyphens]{url}
\usepackage{color}

\usepackage{amsmath,amssymb,amsfonts,amsthm, mathtools, bm}
\usepackage{thmtools, thm-restate}
\usepackage{nicefrac}
\usepackage{array}

\usepackage[normal]{threeparttable}
\usepackage{threeparttablex}
\usepackage{tabularx}
\usepackage{longtable}
\usepackage{booktabs}
\usepackage{colortbl}
\usepackage{multirow}
\usepackage{makecell}

\usepackage[ruled, vlined, linesnumbered]{algorithm2e} 
\SetKw{Continue}{continue}

\usepackage[acronym]{glossaries}

\usepackage{subcaption}
\usepackage{floatrow}
\usepackage{graphicx}
\usepackage{placeins}

\usepackage{tikz}
\usetikzlibrary{external}
\usetikzlibrary{arrows.meta}
\usetikzlibrary{math}
\usetikzlibrary{shapes}
\usepackage{pgfplots}
\usepgfplotslibrary{groupplots}
\usepgfplotslibrary{dateplot}
\pgfplotsset{compat=1.18} 
\usepackage{tikzscale}

\usepackage{enumerate}
\usepackage{multicol}
\usepackage[author=Patrick]{pdfcomment}
\usepackage{xparse}
\usepackage{twoopt}
\usepackage{etoolbox}

\usepackage{eurosym}
\usepackage{ellipsis}
\usepackage{natbib}
\usepackage{paralist}
\usepackage{ulem}
\AtBeginDocument{
	\newgeometry{
		textheight=9in,
		textwidth=6.5in,
		top=1in,
		headheight=14pt,
		headsep=25pt,
		footskip=30pt
	}
}
\newsavebox{\abstractbox}
\renewenvironment{abstract}
{\begin{lrbox}{0}\begin{minipage}{\textwidth}
			\begin{center}\normalfont\sectfont\abstractname\end{center}\quotation}
		{\endquotation\end{minipage}\end{lrbox}%
	\global\setbox\abstractbox=\box0 }

\makeatletter
\expandafter\patchcmd\csname\string\maketitle\endcsname
{\vskip\z@\@plus3fill}
{\vskip\z@\@plus2fill\box\abstractbox\vskip\z@\@plus1fill}
{}{}
\makeatother

\addtokomafont{captionlabel}{\footnotesize\bfseries} 
\addtokomafont{caption}{\footnotesize\bfseries} 

\bibpunct[, ]{(}{)}{,}{a}{}{,}%
\counterwithin*{equation}{section}

\newenvironment{myfont}{\fontfamily{ccr}\selectfont}{\par}
\DeclareTextFontCommand{\textmyfont}{\myfont}

\AtBeginDocument{%
	\abovedisplayskip=7.5pt plus 0pt minus 0pt
	\abovedisplayshortskip=7.5pt plus 0pt minus 0pt
	\belowdisplayskip=7.5pt plus 0pt minus 0pt
	\belowdisplayshortskip=7.5pt plus 0pt minus 0pt
}

\newcolumntype{L}[1]{>{\raggedright\let\newline\\\arraybackslash\hspace{0pt}}p{#1}}
\newcolumntype{C}[1]{>{\centering\let\newline\\\arraybackslash\hspace{0pt}}p{#1}}
\newcolumntype{R}[1]{>{\raggedleft\let\newline\\\arraybackslash\hspace{0pt}}p{#1}}
\renewcommand{\emph}[1]{\textit{#1}}

\newcommand{\emphBold}[1]{{\normalfont \bfseries #1}}
\newcommand{\emphBoldUnderlined}[1]{\underline{{\normalfont \bfseries #1}}}

\begin{document}
\emergencystretch 3em
\newacronym{amod}{AMoD}{autonomous \gls{mod}}
\newacronym{ages}{AGES}{adaptive guided ejection search}


\newacronym{darp}{DARP}{dial-a-ride problem}


\newacronym{ges}{GES}{guided ejection search}


\newacronym{ils}{ILS}{iterated local search}


\newacronym{kdspp}{k-dSPP}{k-disjoint shortest path problem}

\newacronym{ls}{LS}{local search}
\newacronym{lns}{LNS}{large neighborhood search}

\newacronym{mod}{MoD}{mobility-on-demand}


\newacronym{ofat}{OFAT}{one-factor-at-a-time}
\newacronym{osm}{OSM}{OpenStreetMap}

\newacronym{pdptw}{PDPTW}{pickup and delivery problem with time windows}

\newacronym{rnr}{R\&R}{ruin \& recreate}

\newacronym{sa}{SA}{simulated annealing}
\newacronym{sp}{SP}{set partitioning}




\newacronym{wsc}{WSC}{weighted set covering}
\newacronym{wsp}{WSP}{weighted set partitioning}




\newcommand{\mMaxIterations}[1]{M_{#1}}
\newcommand{\mMaxIterationsSplit}{\mMaxIterations{S}}
\newcommand{\mPartialSolutions}{P}
\newcommand{\mVehicleCapacity}{Q}
\newcommand{\mThreshold}{T}
\newcommand{\mThresholdInit}{T^{\textsc{init}}}
\newcommand{\mThresholdDec}{T^{\textsc{dec}}}
\newcommand{\mPerturbationSize}[1]{Z_{#1}}
\newcommand{\mNodeServicetime}[1]{b_{#1}}
\newcommand{\mArcCost}[2]{c_{#1#2}}
\newcommand{\mDeliveryNode}[0]{d}
\newcommand{\mNodeReady}[1]{e_{#1}}
\newcommand{\mNodeDue}[1]{l_{#1}}
\newcommand{\mPickupNode}[0]{p}
\newcommand{\mNodeDemand}[1]{q_{#1}}
\newcommand{\mRequest}[1]{r_{#1}}
\newcommand{\mArcTime}[2]{t_{#1#2}}
\newcommand{\mDirectArcTime}[1]{\mArcTime{\mPickupNode_{#1}}{\mDeliveryNode_{#1}}}
\newcommand{\mVehicleNode}{v}
\newcommand{\mSetArcs}{\mathcal{A}}
\newcommand{\mSetArcsKDSP}{\widetilde{\mathcal{A}}}
\newcommand{\mBlocks}{\mathcal{B}}
\newcommand{\mSetBlocksKDSP}{\widetilde{\mathcal{B}}}
\newcommand{\mSetDeliveries}{\mathcal{D}}
\newcommand{\mSetHyperEdges}{\mathcal{E}}
\newcommand{\mGraph}{\mathcal{G}}
\newcommand{\mGraphKDSP}{\widetilde{\mathcal{G}}}
\newcommand{\mHyperGraph}{\mathcal{H}}
\newcommand{\mSetVehicles}{\mathcal{K}}
\newcommand{\mSetVehiclesKDSP}{\widetilde{\mathcal{K}}}
\newcommand{\mMatching}{\mathcal{M}}
\newcommand{\mSetNodes}{\mathcal{N}}
\newcommand{\mSetNodesKDSP}{\widetilde{\mathcal{N}}}
\newcommand{\bigO}[1]{\ensuremath{\mathcal{O}({#1})}}
\newcommand{\mSetPickups}{P}
\newcommand{\mSetRequests}{R}
\newcommand{\mTimehorizon}{\mathcal{T}}
\newcommand{\mUnassigned}{\ensuremath{\mathcal{U}}}
\newcommand{\mSetHyperNodes}{\mathcal{V}}

\newcommand{\mKDSP}{\ensuremath{\Pi}}
\newcommand{\mTimeLimit}{\Omega_A}
\newcommand{\mProbRelocateExchange}{\ensuremath{\alpha^{\textsc{ils}}}}
\newcommand{\mProbRuinMode}{\ensuremath{\alpha^{\textsc{r\&r}}}}
\newcommand{\mRuinBeta}{\ensuremath{\beta^{\textsc{r\&r}}}}
\newcommand{\mBuffer}{\delta}
\newcommand{\mHyperEdge}[1]{\ensuremath{\varepsilon_{#1}}}
\newcommand{\mArcCostKDSP}[2]{\theta_{#1#2}}
\newcommand{\mRoute}[1]{\ensuremath{\vartheta^{#1}}}
\newcommand{\mRoutePosition}[1]{\ensuremath{\mRoute{}({#1})}}
\newcommand{\mPenaltyDecay}{\lambda}
\newcommand{\mKdspDistanceLimit}{\ensuremath{\mu^{\textsc{dist}}}}
\newcommand{\mKdspTimeLimit}{\ensuremath{\mu^{\textsc{time}}}}
\newcommand{\mVisit}{\nu}
\newcommand{\mUnassignedPenalty}{\ensuremath{\xi}}
\newcommand{\mPermutation}[0]{\ensuremath{\pi}}
\newcommand{\mPermutationPosition}[1]{\ensuremath{\mPermutation({#1})}}
\newcommand{\mPenaltyCounter}{\rho}
\newcommand{\mSolution}{\ensuremath{\sigma}}
\newcommand{\mTime}{\tau}
\newcommand{\mArrivalTime}[1]{\mTime^{\textsc{arr}}({#1})}
\newcommand{\mDepartureTime}[1]{\mTime^{\textsc{dep}}({#1})}
\newcommand{\mAvgSplitSize}{\chi}
\newcommand{\mCardinality}[1]{|{#1}|}
\newcommand{\mSum}[2]{\ensuremath{\sum_{#1}^{#2}}}

\newcommand{\mNodeTimeWindow}[1]{[\mNodeReady{#1},\mNodeDue{#1}]}
\newcommand{\mPred}[1]{\ensuremath{\textsc{Pred}({#1})}}
\newcommand{\mSucc}[1]{\ensuremath{\textsc{Succ}({#1})}}
\newcommand{\mVehicleAssignment}[1]{\ensuremath{\textsc{Assign}({#1})}}
\newcommand{\mLast}[1]{\ensuremath{\textsc{Last}({#1})}}
\newcommand{\mRefFw}[1]{\ensuremath{\textsc{EvalFw}({#1})}}
\newcommand{\mRefBw}[1]{\ensuremath{\textsc{EvalBw}({#1})}}

\newcommand{\mMax}[1]{\ensuremath{\max\{{#1}\}}}
\newcommand{\mMin}[1]{\ensuremath{\min\{{#1}\}}}
\newcommand{\mSequence}{\ensuremath{{\mRoute{}}}}

\newcommand{\mRefCapacitySum}[1]{\ensuremath{q_{\textsc{sum}}({#1})}}
\newcommand{\mRefCapacityMax}[1]{\ensuremath{q_{\textsc{max}}({#1})}}

\newcommand{\mRefTraveltime}[1]{\ensuremath{T_{\textsc{t}}({#1})}}
\newcommand{\mRefEarliestCompletion}[1]{\ensuremath{T_{\textsc{ec}}({#1})}}
\newcommand{\mRefLatestStart}[1]{\ensuremath{T_{\textsc{ls}}({#1})}}
\newcommand{\mRefTWFeasible}[1]{\ensuremath{T_{\textsc{f}}({#1})}}

\newcommand{\mRefCost}[1]{\ensuremath{C_{\textsc{}}({#1})}}

\newcommand{\mMatchingGreedy}{\ensuremath{\textsc{greedy}}}
\newcommand{\mMatchingWSCGreedy}{\ensuremath{\textsc{WSC}^{\textsc{LP}}_{\textsc{greedy}}}}
\newcommand{\mMatchingWSCSRR}{\ensuremath{\textsc{WSC}^{\textsc{LP}}_{\textsc{RR}}}}
\newcommand{\mMatchingWSPGreedy}{\ensuremath{\textsc{WSP}^{\textsc{LP}}_{\textsc{greedy}}}}
\newcommand{\mMatchingWSPSRR}{\ensuremath{\textsc{WSP}^{\textsc{LP}}_{\textsc{RR}}}}

\newcommand{\mPooling}[0]{\ensuremath{p}}
\newcommand{\mPoolings}[0]{\ensuremath{\mathcal{P}}}
\newcommand{\mPoolingCardinality}[0]{\ensuremath{\gamma}}
\newcommand{\mMaxPoolingCardinality}[0]{\ensuremath{\mPoolingCardinality^{\textsc{max}}}}
\newcommand{\mPotentialRequestsForPooling}[1]{\ensuremath{N_{#1}}}
\newcommand{\mMaxInSet}[2]{\underset{#1}{\max}\{#2\}}
\newcommand{\mMinInSet}[2]{\underset{#1}{\min}\{#2\}}

\newcommand{\mSolverSequentialMathHeur}[0]{\ensuremath{\text{MATH}}}
\newcommand{\mSolverDBILS}[0]{\ensuremath{\text{ILS}}}
\newcommand{\mSolverCombined}[0]{\ensuremath{\text{COMB}}}


\title{\large A decomposition-based approach for large-scale pickup and delivery problems}

\author[1]{\normalsize Gerhard Hiermann}
\author[2]{\normalsize Maximilian Schiffer}
\affil{\small 
	TUM School of Management, Technical University of Munich, 80333 Munich, Germany
	
	\scriptsize gerhard.hiermann@tum.de
	
	\small
	\textsuperscript{2}TUM School of Management \& Munich Data Science Institute,
	
	Technical University of Munich, 80333 Munich, Germany
	
	\scriptsize schiffer@tum.de}

\date{}

\lehead{\pagemark}
\rohead{\pagemark}

\begin{abstract}
\begin{singlespace}
{\small\noindent With the advent of self-driving cars, experts envision autonomous mobility-on-demand services in the near future to cope with overloaded transportation systems in cities worldwide. 
Efficient operations are imperative to unlock such a system's maximum improvement potential. Existing approaches either consider a narrow planning horizon or ignore essential characteristics of the underlying problem. In this paper, we develop an algorithmic framework that allows the study of very large-scale pickup and delivery routing problems with more than 20 thousand requests, which arise in the context of integrated request pooling and vehicle-to-request dispatching. We conduct a computational study and present comparative results showing the characteristics of the developed approaches. Furthermore, we apply our algorithm to related benchmark instances from the literature to show the efficacy. 
Finally, we solve very large-scale instances and derive insights on upper-bound improvements regarding fleet sizing and customer delay acceptance from a practical perspective.
\smallskip}
\end{singlespace}
\end{abstract}

\maketitle
\section{Introduction}\label{section:introduction}
Cities worldwide struggle with overloaded transportation systems and related negative externalities: CO\textsubscript{2} emissions cause environmental harm by contributing to the greenhouse effect; health hazards arise from particulate matter and NO\textsubscript{x} emissions; and economic harm stems from congestion-induced lost working hours \citep{Pishue2023}. Over the last decade, the sharing economy paradigm stimulated new mobility services, above all \gls{mod} services to serve individual passenger transportation requests in an urban context. With the advent of self-driving cars, experts envision these \gls{mod} services to be operated by autonomous vehicles in the near future. The resulting \gls{amod} systems promise to enable more efficient on-demand services that are accessible by a large public as they can be offered at a lower cost. Municipalities, practitioners, and scientists share high hopes that such \gls{amod} systems will contribute significantly to reducing the above mentioned negative externalities by allowing, among others, for efficient pooling of passenger requests, congestion aware routing, and convenient feeding to public transport lines \citep{SalazarLanzettiEtAl2019}. Yet, there is so far no consensus whether these benefits will finally lead to reduced negative externalities, or if such a convenient on-demand service may lead to induced demand, such that externalities will be reduced per passenger but not in total \citep{OhSeshadriEtAl2020}.

Independent of this debate, there is consensus that efficient operations are imperative to unlock an \gls{amod} system's maximum improvement potential. To do so, one needs to develop effective, possibly anticipative, control algorithms for the related operational planning tasks: pooling--if possible--passenger requests to shared rides, dispatching these rides to vehicles, and rebalancing vehicles to anticipate future demand. Accordingly, these planning problems sparked the interest of scientists in the field of optimal control, transport optimization, and operations research and led to various algorithmic approaches that range from online control algorithms, often following a (partial) planning-task decomposition scheme \citep{EndersHarrisonEtAl2023, EndersHarrisonEtAl2024}, to integrated algorithms that aim to solve all planning tasks at once in a full information setting to identify an upper bound on the possible system improvement. 

In this context, existing approaches usually suffer from at least one shortcoming: (decomposed) online control algorithms scale to large real-world instances but implicitly assume the performance loss of sequentially deciding on the respective planning tasks to be limited without further discussion \citep{Alonso-MoraSamaranayakeEtAl2017}. Consequently, system improvement analyses are limited to comparisons against the status quo or (naive) baselines. While this setting reveals many interesting insights, it is not sufficient to rigorously analyze the upper bounds of a system's improvement potential, which is of interest for tactical and strategic transportation system analysis and planning. Algorithms that solve the integrated planning problem in a full information setting are often limited in scalability \citep{DoernerSalazar-Gonzaalez2014,SartoriBuriol2020}, such that the obtained full information bounds are not interesting from a practitioner's perspective as the limited problem size diminishes the meaningfulness of the derived results. In the field of operations research, some first attempts exist to solve large-scale pickup and delivery problems in a full information setting. However, these works are motivated by freight transport applications, such that the studied instance size still remains below the scale encountered in an \gls{amod} context. Moreover, the instances studied in this work reveal significantly different characteristics with respect to demand and time windows such that it remains questionable whether existing--usually highly-tailored--algorithms will provide a good solution quality within an \gls{amod} context.

Against this background, we aim at developing an algorithmic framework that allows to solve very large-scale pickup and delivery routing problems that arise in the context of integrated request pooling and vehicle-to-request dispatching in \gls{amod} systems. Focusing on offline full information problem settings, we design this algorithmic framework in such a way that it allows us to analyze not only a full information bound of the system's performance improvement but also the impact of taking pooling and dispatching decisions sequentially or in an integrated fashion.

In the following, we first briefly review related literature before we specify our contribution and the paper's organization.

\subsection{State of the Art}

From a general perspective, controlling a ride-hailing fleet is related to solving an \gls{pdptw} \citep{SavelsberghSol1995, RopkeCordeauEtAl2007}. Literature on the \gls{pdptw} focuses overwhelmingly on deliveries of goods \citep{BattarraCordeauEtAl2014}. Apart from its basic problem variant, several extensions for the \gls{pdptw} have been studied, e.g., focusing on loading constraints \citep{IoriMartello2010}, selection of requests \citep{Al-ChamiManierEtAl2016}, or recharging of electric vehicles \citep{Goeke2019}. Passenger transportation variants are studied under the umbrella of the \gls{darp}, which was conceived to formalize a mobility service for elderly and disabled people to maximize the quality of service while minimizing cost. Herein, not only punctuality in service but the time passengers travel in the vehicle are considered using the aptly called user-ride-time constraint to limit the ride duration and avoid long periods of waiting times in the vehicle \citep{DoernerSalazar-Gonzaalez2014}. 
This work focuses on inner-city ride-sharing operations with tight time windows regarding pickup and dropoff and short trips. As such, the maximum ride times of users are bounded by the time windows, removing the necessity of using explicit limits. For this reason, we consider the \gls{pdptw} as our modeling basis and focus our literature review on this problem. 

Only a few exact approaches exist to solve the \gls{pdptw}. Besides some early branch-and-price solutions \citep{DumasDesrosiersEtAl1991, SavelsberghSol1998}, branch-and-price-and-cut algorithms have been developed by \cite{RopkeCordeau2009} and \cite{BaldacciBartoliniEtAl2011} with \cite{BettinelliCeselliEtAl2014} and \cite{GschwindIrnichEtAl2018} focusing on how to apply bidirectional search in the pricing problem. The largest instances that could be solved with these approaches comprise 100 requests for regular instances and up to 500 requests for tightly constrained instances. Recently, \cite{VadsethAnderssonEtAl2023} proposed a route modifying improvement model, which they initialized with best-known solutions to find improvements for the large benchmark sets.

To solve larger instances, different metaheuristic resolution approaches have been developed, covering Tabu Search \citep{NanryBarnes2000, LiLim2003}, genetic algorithm \cite{Pankratz2005} and hyper-heuristics \citep{NasiriKeedwellEtAl2022}. \cite{NagataKobayashi2010} proposed a \gls{ges} and focused on minimizing the fleet only.
Current state-of-the-art methods include an iterative approach by \cite{CurtoisLandaSilvaEtAl2018}, using an \gls{ages}, \gls{lns} and \gls{ls} sequentially to repeatedly probe minimizing the fleet before optimizing costs. \cite{SartoriBuriol2020} later improved upon this approach and proposed a matheuristic composing of a \gls{ages}, \gls{lns}, and an \gls{sp} component. The latter continuously tries to recombine routes of solutions found during the search to yield better solutions. These components are nested in an \gls{ils} framework to diversify the search.
\cite{ChristiaensVandenBerghe2020} proposed an \gls{rnr} using a novel slack-induced string removal and greedy insertion with blinks operator. The approach showed excellent performance on various problems, including the \gls{pdptw}.

In contrast with the literature, our focus lies on very large instances, ranging up to 21 thousand requests compared to the currently studied 2500 available requests.

\subsection{Contribution}
With this work, we provide a new state of the art for solving very large scale \glspl{pdptw} in the context of ride-sharing. We provide an algorithmic framework that comprises a decomposition-based matheuristic which allows to solve instances with up to 5000 requests in a few minutes, an \gls{ils}-based metaheuristic that yields a better solution quality within a computational time limit of up to 15 minutes, as well as a hybrid approach that uses our matheuristic to warm-start the \gls{ils} in order to further improve solution quality. 

With this algorithmic framework, we provide a thorough computational study to compare the proposed algorithms against each other and understand the respective algorithmic characteristics. We then use our metaheuristic algorithm to solve the benchmark data set for the \gls{pdptw} in the context of ridesharing \cite{SartoriBuriol2020}. Here, we show that our algorithm improves significantly over the algorithm of \cite{SartoriBuriol2020} and find new best-known solutions. Lastly, we apply our algorithm to study very large-scale instances with up to 21375 requests that have not been solved before. By so doing, we shed light on upper-bound improvements with respect to fleet sizing and customer delay acceptance from a practical perspective.

\subsection{Organization}
The remainder of this paper is organized as follows: Section~\ref{section:problem_formulation} formally introduces our problem setting. Section~\ref{section:methodology} details the proposed solution approach. In Section~\ref{section:experimental-design}, we outline our experimental design and present the related numerical results in Section~\ref{section:results}. Finally, Section~\ref{section:conclusion} concludes this work with a short summary and an outlook on future research.

\section{Problem Setting}\label{section:problem_formulation}
We study an offline problem setting in which a fleet operator has full knowledge about requests that arrive in a certain time horizon, e.g., a peak hour, an operational shift, or a day. Formally, we denote this time horizon by $\mTimehorizon$. The fleet operator operates a fleet of constant size and offers a ride-hailing service. In this context the operator decides on i) which customer requests to pool to a shared ride, and ii) which vehicle to dispatch for operating each (shared) ride. Clearly, the operator can also decide to not dispatch a vehicle to a ride, which implicitly models request rejections.

The fleet operates on a road network and we model respective operations on a fully connected graph $\mGraph=(\mSetNodes,\mSetArcs)$ composed of a set of nodes $\mSetNodes = \mSetPickups \cup \mSetDeliveries \cup \mSetVehicles$ and a set of arcs $\mSetArcs$. Here, Set $\mSetPickups$ contains the pickup nodes of all passenger requests, $\mSetDeliveries$ contains the respective drop-off nodes, and set $\mSetVehicles$ contains vehicle nodes, i.e., nodes that indicate the initial positions of all vehicles. Arcs $(i,j) \in \mSetArcs$ represent paths through the street network with corresponding cost $\mArcCost{i}{j}$ and travel time~$\mArcTime{i}{j}$.

We represent a passenger request as a quadruple $(\mPickupNode_r,\mDeliveryNode_r,\mNodeReady{r},\mNodeDue{r})$, where $\mPickupNode_r \in \mSetPickups$ and $\mDeliveryNode_r \in \mSetDeliveries$ denote the pickup and dropoff location of request $r\in R$, and $\mNodeReady{r},\mNodeDue{r}$ define the request's time window $\mNodeTimeWindow{r}$. This time window indicates the earliest time $\mNodeReady{r} \in \mTimehorizon$ at which it is possible to pick up the passenger at $\mPickupNode_r$ and the latest time $\mNodeDue{r} \in \mTimehorizon$ at which the passenger needs to be dropped at $\mDeliveryNode_r$. In this context, one can interpret $\mNodeReady{r}$ either as the time at which a passenger sends a request to the operator, or as a specified pickup time that lies further in the future. To account for a passengers willingness to accept a detour when participating in a shared ride, e.g., incentivized by a price reduction, we calculate $\mNodeDue{r}$ as
$$\mNodeDue{r} = \mNodeReady{r} + \mArcTime{\mPickupNode_r}{\mDeliveryNode_r} + \mBuffer_r,$$
where $\mArcTime{\mPickupNode_r}{\mDeliveryNode_r}$ is the travel time between $\mPickupNode_r$ and $\mDeliveryNode_r$, and $\mBuffer$ denotes a maximum time budget available for detours. In the remainder of this paper, we will refer to this time budget as a request's buffer.

With this notation, one can easily link a request's time window information to its origin and destination. We do so by defining a pickup time window $\mNodeTimeWindow{\mPickupNode_r}$ as $\mNodeReady{\mPickupNode_r} = \mNodeReady{r}$ and $\mNodeDue{\mPickupNode_r} = \mNodeReady{r} + \mBuffer$, and a delivery time window $\mNodeTimeWindow{\mDeliveryNode_r}$ as $\mNodeReady{\mDeliveryNode_r} = \mNodeReady{r} + \mArcTime{\mPickupNode_r}{\mDeliveryNode_r}$ and $\mNodeDue{\mDeliveryNode_r} = \mNodeDue{r} = \mNodeReady{r} + \mArcTime{\mPickupNode_r}{\mDeliveryNode_r} + \mBuffer$ for each request. While these definitions appear to be redundant in the problem description, they will ease notation and clarity when discussing our algorithmic framework.

\textit{Solution representation:} A solution $\mSolution$ represents a set of routes $\mSolution=\{\mRoute{1},\ldots,\mRoute{|\mSetVehicles|}\}$, one for each vehicle. 
Each route $\mRoute{k}=(\mVisit_0, \ldots, \mVisit_n)$ is a sequence of nodes $\mVisit_i \in \mSetNodes$ that denotes in which order a vehicle visits them. Here, $\mVisit_0$ is always the starting location of the corresponding vehicle. Each sequence implicitly defines departure and arrival times, which can be trivially calculated by propagating travel times, departing as early as possible from each node. We define $\mArrivalTime{\mRoute{},\mVisit_i}$ and $\mDepartureTime{\mRoute{},\mVisit_i}$ to access the arrival and departure time of node $\mVisit_i \in \mSetNodes$. Finally, we define $\mUnassigned(\mSolution)$ as the set of unassigned requests that are not served by any vehicle within a solution $\mSolution$.

\textit{Constraints:} A valid solution $\mSolution$ has to adhere to the following constraints.
\begin{enumerate}[i]
    \item Vehicles can never exceed their capacity $\mVehicleCapacity$ and can thus operate at maximum $\mVehicleCapacity$ requests in parallel.
    \item If served, request nodes $\mPickupNode_{\mRequest{}}$ and $\mDeliveryNode_{\mRequest{}}$ have to be visited by the same vehicle, i.e., \begin{align}
	    \mPickupNode_{\mRequest{}} \in \mRoute{} \leftrightarrow  \mDeliveryNode_{\mRequest{}} \in \mRoute{}, \quad \forall \mRequest{} \in \mSetRequests\label{eq:problem:same-vehicle}
	\end{align}
    \item Request nodes must be served in order and inside the corresponding requests' time window [$\mNodeReady{r}, \mNodeDue{r}$]. If a vehicle arrives early, it has to wait until the request can be served. Formally \begin{align}
	    \mNodeReady{r} \leq \mDepartureTime{\mPickupNode_r} \leq \mArrivalTime{\mDeliveryNode_r} \leq \mNodeDue{r},  \quad \forall \mRequest{} \in \mSetRequests,\label{eq:problem:precedence}
	\end{align}
\end{enumerate}

\textit{Objective function:} We consider a hierarchical objective, minimizing two quantities: i) the number of unserved requests in $\mUnassigned(\mSolution)$, and ii) the total travel cost considering the driving distance $\mArcCost{i}{j}$.
\begin{align}
\begin{split}
\min & \quad |\mUnassigned(\mSolution)|
\\
\min & \quad \sum_{\mRoute{} \in \mSolution}^{} \sum_{i=1}^{|\mRoute{}|} \mArcCost{\mVisit_{i-1}}{\mVisit_{i}} 
\label{eq:problem:objective}
\end{split}
\end{align} 

Among all feasible solutions fulfilling these constraints, we seek a solution $\mSolution^*$ that minimizes the objective function \eqref{eq:problem:objective}.

\textit{Discussion:} Three comments on our problem setting are in order. First, we limit our problem setting to a full information scenario, which omits directly leveraging our algorithm for fleet control in practice. While limiting, this simplification is in line with our paper's scope: solving large-scale full information instances to obtain upper bounds on the system's improvement potential, as well as a temporally unbiased analyses on whether decomposing pooling and dispatching decisions affects solution quality or not. Second, while we analyze the impact of decoupled or integrated pooling and dispatching decisions, we ignore explicit rebalancing and rerouting. Omitting rebalancing is reasonable as it is only beneficial during online decision making but not in a full information setting. For (congestion aware) rerouting, recent works show significant improvement potential even if it is conducted on a subsequent decision level (Jalota et al 2023). Accordingly, we ignore this aspect to isolate the effect of decomposing or integrating pooling and dispatching decisions. Lastly, we like to mention that the proposed problem setting differs from classical \glspl{pdptw} in three aspects: i) we do not consider heterogeneous demand, instead each request has a demand of one, i.e., represents one customer; ii) our setting resembles an orienteering problem where vehicles are initially scattered within the service area instead of being located at a central depot and may also end their last service at an arbitrary location; iii) we consider a different objective function that aims at minimizing the number of unserved requests and subsequently the respective operational costs for a given fleet size. All of these differences result from our ride-hailing application, which significantly differs from the usually studied logistics context.

\section{Methodology}\label{section:methodology}
This section details our algorithmic framework. To obtain a framework that allows to take decomposed as well as integrated pooling and dispatching decisions, we proceed as follows. In a first step, we develop a matheuristic that decomposes the planning problem and takes pooling and dispatching decisions sequentially. We then focus on integrated decision-making and develop a metaheuristic that allows to take integrated pooling and dispatching decisions.

\subsection{Sequential Pooling \& Dispatching}\label{section:methodology:matheuristic}
\begin{figure}[b]
	\centering
\begingroup%
  \makeatletter%
  \providecommand\color[2][]{%
    \errmessage{(Inkscape) Color is used for the text in Inkscape, but the package 'color.sty' is not loaded}%
    \renewcommand\color[2][]{}%
  }%
  \providecommand\transparent[1]{%
    \errmessage{(Inkscape) Transparency is used (non-zero) for the text in Inkscape, but the package 'transparent.sty' is not loaded}%
    \renewcommand\transparent[1]{}%
  }%
  \providecommand\rotatebox[2]{#2}%
  \newcommand*\fsize{\dimexpr\f@size pt\relax}%
  \newcommand*\lineheight[1]{\fontsize{\fsize}{#1\fsize}\selectfont}%
  \ifx\svgwidth\undefined%
    \setlength{\unitlength}{327.12253096bp}%
    \ifx\svgscale\undefined%
      \relax%
    \else%
      \setlength{\unitlength}{\unitlength * \real{\svgscale}}%
    \fi%
  \else%
    \setlength{\unitlength}{\svgwidth}%
  \fi%
  \global\let\svgwidth\undefined%
  \global\let\svgscale\undefined%
  \makeatother%
  \begin{picture}(1,0.27948586)%
    \lineheight{1}%
    \setlength\tabcolsep{0pt}%
    \put(0,0){\includegraphics[width=\unitlength,page=1]{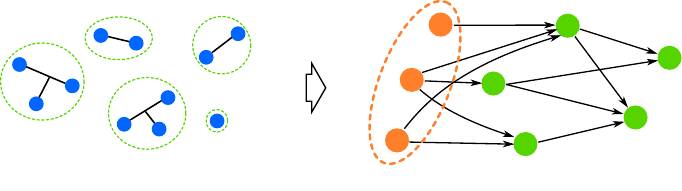}}%
    \put(0.51462977,0.00093328){{\makebox(0,0)[lt]{\lineheight{1.25}\smash{\begin{tabular}[t]{l}vehicles\end{tabular}}}}}%
    \put(0.71761427,0.00093328){{\makebox(0,0)[lt]{\lineheight{1.25}\smash{\begin{tabular}[t]{l}matched requests\end{tabular}}}}}%
  \end{picture}%
\endgroup%

	\caption{Outline of the sequential pooling and dispatching approach.}
	\begin{subfigure}[c]{0.05\textwidth}\end{subfigure}
	\begin{subfigure}{0.40\textwidth}
		\subcaption{Step 1: Hypergraph matching.}
		\label{fig:sequential-pooling-and-dispatching:pooling}
	\end{subfigure}
	\begin{subfigure}[c]{0.05\textwidth}\end{subfigure}
	\begin{subfigure}{0.46\textwidth}
		\subcaption{Step 2: Create and solve a dispatching graph.}
		\label{fig:sequential-pooling-and-dispatching:dispatching}
	\end{subfigure}
	\label{fig:sequential-pooling-and-dispatching}
\end{figure}

To devise an algorithm for sequential pooling and dispatching decisions, we expand two of our recent works that focused on the respective isolated decision tasks. Figure~\ref{fig:sequential-pooling-and-dispatching} illustrates the rationale of our algorithm. 
In a first step (Figure~\ref{fig:sequential-pooling-and-dispatching:pooling}), we model potential request poolings as hyperedges in a hypergraph, which allows us to effectively pool customer requests by obtaining a maximum weighted matching. We then model these pooled requests in a dispatching graph (Figure~\ref{fig:sequential-pooling-and-dispatching:dispatching}), which allows us to calculate the request to vehicle assignments in polynomial time.

\subsubsection{Generating feasible hyperedges}

We define a hypergraph $\mHyperGraph = (\mSetHyperNodes, \mSetHyperEdges)$ with a vertex set $\mSetHyperNodes$, a hyperedge set $\mSetHyperEdges$, and its hyperedge weights $\omega$, which relates to our problem as follows:
\begin{enumerate}[i]
\item The vertex set is ordered $\mSetHyperNodes = [n]$, (with $n = |\mSetHyperNodes|$) and each vertex represents a request $\mRequest{} \in \mSetRequests$. The set is sorted by the earliest service time $\mNodeReady{r}$ in ascending order. 
\item Each hyperedge $\mHyperEdge{}$ with $|\mHyperEdge{}| \geq 2$ represents a potential pooled ride that contains multiple requests.
\item The edge weights $\omega: \mSetHyperEdges \to \mathbb{R}_+$, represent the utility of a potentially pooled ride.
\end{enumerate}
The number of hyperedges can be potentially intractable, with up to $\mSum{\gamma=2}{|\mSetRequests|} \binom{|\mSetRequests|}{\gamma}$ possible combinations. To tackle this issue, we limit the hyperedges generated in two ways. First, we limit the rank of the hypergraph to $\mMaxPoolingCardinality = 4$, i.e., $|\mHyperEdge{}| \leq \mMaxPoolingCardinality, \forall \mHyperEdge{} \in \mSetHyperEdges$. Second, we define $\mPotentialRequestsForPooling{r}$ as the set of feasible neighbors of request $r$, where $\mNodeReady{r} \leq \mNodeReady{r'} \leq \mNodeDue{r} + \mBuffer$. For each request $r$ and its neighbors $N_r$, we generate every combination $\pi \in \Pi_r$ of up to $\mMaxPoolingCardinality$ requests, $\pi \in \{r\}\cup \mPotentialRequestsForPooling{r}$. We only allow $r' \in \mPotentialRequestsForPooling{r}$ if $r'$ is positioned after $r$ in the ordered vertex set $\mSetHyperNodes$ to avoid symmetries.

Note that not all possible hyperedge combinations may result in a feasible sequence as they may violate the pickup time feasibility or capacity constraint, and can be ignored. Furthermore, a set of requests can be sequenced in various ways while respecting the precedence constraint. For example, $\pi = \{r, r'\}$ can be sequenced in four ways: $(\mPickupNode_{r}, \mPickupNode_{r'}, \mDeliveryNode_{r}, \mDeliveryNode_{r'})$, $(\mPickupNode_{r}, \mPickupNode_{r'}, \mDeliveryNode_{r'}, \mDeliveryNode_{r})$, $(\mPickupNode_{r'}, \mPickupNode_{r}, \mDeliveryNode_{r}, \mDeliveryNode_{r'})$, $(\mPickupNode_{r'}, \mPickupNode_{r}, \mDeliveryNode_{r'}, \mDeliveryNode_{r})$. 
To limit the computational complexity, we reduce the number of possible sequences for each $\pi \in \Pi_{r}^{\gamma}$ by only considering the cheapest feasible sequence of visits in terms of travel cost.

\subsubsection{Hypergraph matching}
\label{section:methodology:matheuristic:matching}

We define a matching $\mMatching$ for hypergraph $\mHyperGraph$ as a subset of hyperedges $\mMatching \subset \mSetHyperEdges$ where all edges are disjoint. Such a matching is maximal if it is not a strict subset of any other matching, and it is maximum if no other matching with a greater cardinality exists. Then, the best pooling strategy $\mMatching^*$ equals a maximum weight matching in $\mHyperGraph$, formally \begin{align}
    &\mMatching^*:= \text{arg}~\underset{X \in \mSetHyperEdges}{\text{max}} \mSum{\mHyperEdge{} \in \mSetHyperEdges}{} \omega(\mHyperEdge{})\nonumber\\
    &\quad \text{s.t.} \quad \mHyperEdge{i} \cup \mHyperEdge{j} = \emptyset, \forall \mHyperEdge{i}, \mHyperEdge{j} \in \mSetHyperEdges \nonumber
\end{align}

Clearly, the definition of the hyperedge weights $\omega$ is crucial for the actual matching performance. In this work, we analyze four weight functions to evaluate hyperedges $\mHyperEdge{} \in \mSetHyperEdges$ as follows.
\begin{align}
\omega^1(\mHyperEdge{}) &= -\mMaxPoolingCardinality / \mPoolingCardinality_{\mHyperEdge{}}\\
\omega^2(\mHyperEdge{}) &= (\sum_{(i,j) \in \mHyperEdge{}}^{} \mArcCost{i}{j}) \cdot \omega^1\\
\omega^3(\mHyperEdge{}) &= ((\mMaxInSet{r \in \mHyperEdge{}}{\mNodeDue{r}}-\mMinInSet{r \in \mHyperEdge{}}{\mNodeReady{r}}) - (\mMinInSet{r \in \mHyperEdge{}}{\mNodeDue{r}}-\mMaxInSet{r \in \mHyperEdge{}}{\mNodeReady{r}})) \cdot \omega^1\\
\omega^4(\mHyperEdge{}) &= \omega^2 (1-\rho) + \omega^3 \rho
\end{align}

Here, $\omega^1(\mHyperEdge{})$ parameterizes each edge based on the negative inverse of its cardinality normalized by the maximum cardinality of hyperedge $\mHyperEdge{}$. While this negative inverse seems unintuitive at first sight, it appears reasonable once we embedd the respective poolings in our dispatching algorithm, which bases on solving a \gls{kdspp} on a dispatching graph. To do so, it requires a negative representation of each poolings benefit. We reach such a notion with respect to each poolings utilization by considering the negative inverse of a utilization ratio expressed via the hyperedge cardinality. In the remaining weight definitions, we use $\omega_1$ as a normalization factor such that no further transformations are necessary to define the other weights in a metric that is suitable for the subsequent dispatching algorithm. The remaining weights aim to either account for the cost of the resepctive request sequence ($\omega^2$), the temporal overlap between the pooled requests ($\omega^3$), or a convex combination of both ($\omega^4$).

To find $\mMatching$, we use a two-step matheuristic approach. First, we solve a continuous \gls{wsc} problem to select a subset of hyperedges that may be connected to the same request nodes, i.e., a fractional matching $\overline{\mMatching}$. We then derive an integral matching $\mMatching$ based on the solution of the \gls{wsc} using a greedy selection procedure. 

We define our continuous \gls{wsc} as follows.
\begin{align}
    \max \quad &\sum_{\mHyperEdge{} \in \mSetHyperEdges'} \omega(\mHyperEdge{}) \cdot x_{\mHyperEdge{}} &&\\
    \text{s.t.} \quad & \nonumber \\
     &\sum_{\mHyperEdge{} \in \mSetHyperEdges'} a_{\mRequest{}\mHyperEdge{}} x_{\mHyperEdge{}} \geq 1 &&\forall \mRequest{} \in \mSetRequests \label{eq:cover}\\
     &0 \leq x_{\mHyperEdge{}} \leq 1 &&\forall \mHyperEdge{} \in \mSetHyperEdges'
\end{align}
Here, $a_{\mRequest{}\mHyperEdge{}} = 1$ if request $\mRequest{}$ is connected to hyperedge $\mHyperEdge{}$ and remains zero otherwise; while $\mSetHyperEdges' = \mSetHyperEdges \cup \bigcup_{r \in \mSetRequests} \mHyperEdge{r}$ is the set of hyperedges including hyperedges comprising only a single request $r \in \mSetRequests$. 

After solving the continuous \gls{wsc} problem, we obtain matching $\mMatching$ from the fractional matching $\overline{\mMatching}$ using a greedy selection procedure.
Herein, we traverse the hyperedges $\mHyperEdge{} \in \overline{\mMatching}$, sorted by their fractional solution value in descending order, to construct a matching $\mMatching$. We only add a hyperedge to our matching if it does not contain any request already covered in $\mMatching$. The process stops after all requests have been selected or the list has been traversed.

We tested additional approaches to find good matchings in preliminary experiments, including a simple greedy selection approach based on the hyperedge weight $\omega(\mHyperEdge{})$, and a weighted set partitioning reformulation (see Appendix~~\ref{appendix:pooling}). Although the aggregated results do slightly favor the greedy approach, the disaggregated results on fleet sizes relevant in practice show better performance when using the matheuristic approach with $\omega^4$ and $\rho=0.7$. 

\subsubsection{Dispatching}
\label{section:methodology:matheuristic:dispatching}

Preliminarily analyses of the characteristics of our planning problem showed that requests in our ride hailing content often tend to interlace less compared to classical \gls{darp} instances. Accordingly, we often observe subsequences of visits that start and end with an empty vehicle, often referred to as zero-split sequences \cite{ParraghDoernerEtAl2010}, fragments \cite{RistForbes2021}, or zero-sum blocks \cite{MalheirosRamalhoEtAl2021}. We design our sequential approach and the respective dispatching algorithm around this property: we create respective zero-sum blocks via the hypergraph matching described in Section~\ref{section:methodology:matheuristic:matching} and use these to construct a dispatching graph that allows us to obtain the respective vehicle-to-request dispatching by solving a \gls{kdspp}. The remainder of this section briefly outlines how we create the respective dispatching graph and a polynomial \gls{kdspp} algorithm. 

\begin{figure}[t]
	\centering
	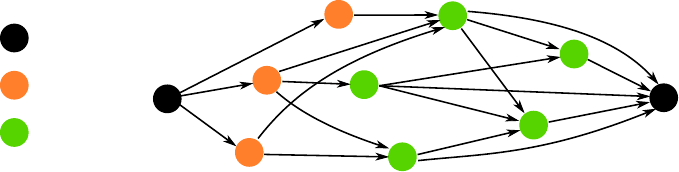
	\caption{Dispatching graph transformed into a k-disjoint shortest path problem}
	\label{fig:dispatching}
\end{figure}

We start by creating our dispatching graph as a weighted and directed source-sink graph $\mGraphKDSP= (\mSetNodesKDSP,\mSetArcsKDSP, \mArcCostKDSP{}{})$, consisting of a set of vertices $\mSetNodesKDSP$, a set of arcs $\mSetArcsKDSP$, and a vector of weights $\mArcCostKDSP{}{}$, containing a weight $\mArcCostKDSP{a}{}$ for each arc $a \in \mSetArcsKDSP$ (see Figure~\ref{fig:dispatching}).
The vertex set $\mSetNodesKDSP = \{w,w^*\} \cup \mSetVehiclesKDSP \cup \mSetBlocksKDSP$, comprises a dummy source $w$ and a dummy sink node $w^*$, a subset of vehicle vertices $\mSetVehiclesKDSP$, and a subset of pooled requests $\mSetBlocksKDSP$ which we refer as blocks. Here, we associate each vehicle vertex with a vehicle's starting position while associating each block vertex with a feasible sequence of pickup and delivery node visits.
More formally, a block $b \in \mSetBlocksKDSP$ is a tuple $(\widetilde{\mSequence_b}, \widetilde{\mNodeReady{}}_{b}, \widetilde{t}_b)$ comprise of a sequence of visits $\widetilde{\mSequence}$ for a set of requests, such that precedence, capacity and time window constraints are satisfied. Next, $\widetilde{\mNodeDue{}}_{b}$ defines the fixed starting time and $\widetilde{t}_b$ the duration to serve $\widetilde{\mSequence}$.

Arcs are constructed as follows: we create an arc between two block vertices if a block-to-block connection does not violate the time window constraint. Similarly, we create an arc between a vehicle and a block vertex if the block can be reached in time by the vehicle. Finally, we connect the dummy source with each vehicle vertex and connect all vehicle and block vertices to the dummy sink. 

A path in $\mGraphKDSP$ that starts at dummy source $w$ and ends in dummy sink node $w^*$ represents a feasible vehicle trip. We can directly associate the trip to the vehicle dependent on the vehicle vertex it contains. We conclude the construction by assigning a weight $\mArcCostKDSP{u}{v}$ to each $(u,v) \in \mSetArcsKDSP$. Herein, we set the arcs that leave the source or enter the sink to zero. Next, we set the weights of vehicle-to-block arcs $(v,\widetilde{b})$ to $\mArcCostKDSP{v}{b} = \mArcCost{v}{\widetilde{\mSequence}_b} - \mUnassignedPenalty \frac{\mCardinality{\widetilde{\mSequence_b}}}{2}$. The second part is the cost traveling from vehicle vertex $v$ to the first visit in $b$; the second part represents the profit of serving the requests of a block, with $\mUnassignedPenalty$ being a large profit term, such that serving requests is always preferred.
The weights of the remaining block-to-block arcs connecting block $b$ to $b'$ are set to $\mArcCostKDSP{b}{b'} = \mArcCost{\widetilde{\mSequence}_b}{\widetilde{\mSequence}_{b'}} - \mUnassignedPenalty \frac{\mCardinality{\widetilde{\mSequence_{b'}}}}{2}$. The first part is the cost of traveling from the last visit in $\widetilde{\mSequence}_{b}$ to the first visit in $\widetilde{\mSequence}_{b'}$, whereas the second part represents again the profit of serving the requests. 

With graph $\mGraphKDSP$ defined, we now briefly discuss the polynomial-time algorithm for the \gls{kdspp}. As proven by \cite{Suurballe1974}, it is possible to increase the number of disjoint paths on graph $\mGraphKDSP$ from $i$ to $i+1$ by finding the shortest interlacing on a modified graph $\mGraphKDSP'$. Herein, we initialize the algorithm by finding the shortest path using the Bellman-Ford algorithm. Note that $\mArcCostKDSP{u}{v}$ could be negative, depending on $\mUnassignedPenalty$. Therefore, we reweight the arcs similar to the well-known Johnson algorithm, using the result from the Bellman-Ford algorithm. Afterward, we continue with a modified Dijkstra algorithm to iteratively find the i+1-disjoint shortest paths until we reach $k=\mCardinality{\mSetVehiclesKDSP}$. For a detailed description of the algorithm, we refer to \cite{Suurballe1974} for the original implementation, or the more recent use in \cite{SchifferHiermannEtAl2021}.

After identifying a set of disjoint shortest paths $\Pi$, we generate a complete solution by concatenating the visit sequence in blocks based on the paths and assigning them to the corresponding vehicle represented by vehicle node $u_c \in \mSetVehiclesKDSP$. Blocks that are not contained in any path remain unassigned.

\subsection{Integrated approach}\label{section:methodology:metaheuristic}
In this section, we introduce an \gls{ils}-based algorithm that allows to take integrated pooling and dispatching decisions. When solving large-scale instances, it is imperative to chose an efficient solution represantation and respective evaluations to obtain a performant algorithm. Accordingly, we first focus in Section~\ref{section:methodology:eval} on the solution representation and evaluation methods used in our algorithm. We then Focus on the \gls{ils}-based algorithm in Section~\ref{section:methodology:ils}, before detailing a ruin and recreate procedure and further intensification techniques in Section~ \ref{section:methodology:rnr}. To keep this paper concise, we focus all discussions on the cost-minimizing objective and detail the respective fleet minimization component in Appendix~\ref{appendix:fleet-min-algorithm}.

\subsubsection{Solution representation and evaluation}\label{section:methodology:eval}

Recall that we represent a solution for our problem as a set of $\left|\mSetVehicles\right|$ routes, one for each vehicle. Typically, one stores each route as a list of nodes that denotes the order in which they are visited. Although the number of required lists is known, the length of each list can vary and depends on each instance's properties, e.g., the number of requests or time horizon, and only materializes throughout the search. In the most extreme case, a single route might contain all requests, thus requiring \bigO{2\mCardinality{\mSetRequests}} memory, such that we obtain an instance-dependent upper bound on a solution's memory requirement by \bigO{\mCardinality{\mSetVehicles}\cdot2\mCardinality{\mSetRequests}}. 

To maintain each list, we exploit that any node is not visited more than once and maintain the following information for each node $i \in \mSetNodes$: while $\mPred{i}$ denotes the predecessor of $i$, $\mSucc{i}$ denotes its successor; $\mVehicleAssignment{i}$ yields the vehicle id node $i$ is assigned to (receiving label $\mCardinality{\mSetVehicles}+1$ if unassigned); $\mRefFw{i}$ holds evaluation information up until node $i$, starting from the vehicle node, whereas $\mRefBw{i}$ holds evaluation information starting from node $i$ to the last visit.

The data in $\mRefFw{i}$ contains the evaluation information of sequence $\{\mVehicleNode{}, \ldots, i\}$ and $\mRefBw{i}$ holds $\{i, \ldots, \mRoutePosition{-1}\}$, repectively, where $\mRoutePosition{-1}$ is the last visit of route $\mRoute{}$. We calculate both iteratively using the distance and time window-related formulas of \cite{VidalCrainicEtAl2012} and demand calculations analogous to \cite{BulhoesSubramanian2018} as follows.

\textit{Capacity evaluation:} Let $q_{sum}(\mSequence)$ be the current number of passengers served, and let $q_{max}(\mSequence)$ be the maximum number of concurrent passengers served in $\sigma$. Given two subsequences $\mSequence^1$ and $\mSequence^2$, we can calculate the concatenated sequence using operation $\oplus$ as follows:
\begin{align}
    \mRefCapacitySum{\mSequence^1 \oplus \mSequence^2} &= \mRefCapacitySum{\mSequence^1} + \mRefCapacitySum{\mSequence^2}\\
    \mRefCapacityMax{\mSequence^1 \oplus \mSequence^2} &= \mMax{\mRefCapacityMax{\mSequence^1},\mRefCapacityMax{\mSequence^1 + \mSequence^2}}
\end{align}
where a sequence $\mSequence'$ containing a single node $i$ is initialized as follows:
\begin{align}
    \mRefCapacitySum{\mSequence'} &= \begin{cases}
1 & \text{ if } i \in \mSetPickups \\
-1 & \text{ if } i \in \mSetDeliveries \\
0 & \text{ if } i \in \mSetVehicles
\end{cases}\\
    \mRefCapacityMax{\mSequence'} &= \mMax{0, \mRefCapacitySum{\mSequence'}}
\end{align}
A route \mRoute{} is feasible regarding the capacity constraint if $\mRefCapacityMax{\mRoute{}} \leq \mVehicleCapacity$. Note that the load of a vehicle cannot be negative as long as its route satisfies the respective the precedence constraint \eqref{eq:problem:precedence}.

\textit{Travel time evaluation:} Let $\mRefTraveltime{\mSequence}$ be the accumulated travel time without potential waiting time, and let $\mRefLatestStart{\mSequence}$, $\mRefEarliestCompletion{\mSequence}$ be the latest start and earliest completion time, respectively; $\mRefTWFeasible{\mSequence}$ indicates whether the sequence is time windows feasible. Then, we define our concatenation operation $\oplus$ as follows:
\begin{align}
    \mRefTraveltime{\mSequence^1 \oplus \mSequence^2} &= \mRefTraveltime{\mSequence^1} + \mArcTime{\mSequence^1}{\mSequence^2} + \mRefTraveltime{\mSequence^2}\\
    \mRefEarliestCompletion{\mSequence^1 \oplus \mSequence^2} &= \mMax{\mRefEarliestCompletion{\mSequence^1} + \mArcTime{\mSequence^1}{\mSequence^2} + \mRefTraveltime{\mSequence^2}, \mRefEarliestCompletion{\mSequence^2}}\\
    \mRefLatestStart{\mSequence^1 \oplus \mSequence^2} &= \mMin{\mRefLatestStart{\mSequence^1}, \mRefLatestStart{\mSequence^2} - \mArcTime{\mSequence^1}{\mSequence^2} - \mRefTraveltime{\mSequence^1}}\\
    \mRefTWFeasible{\mSequence^1 \oplus \mSequence^2} &= \mRefTWFeasible{\mSequence^1} \wedge \mRefTWFeasible{\mSequence^2} \wedge (\mRefEarliestCompletion{\mSequence^1} + \mArcTime{\mSequence^1}{\mSequence^2} \leq \mRefLatestStart{\mSequence^2})
\end{align}

A single-node sequence $\mSequence'$ is initialized with $\mRefTraveltime{\mSequence'} = 0, \mRefEarliestCompletion{\mSequence'} = \mNodeReady{i}, \mRefLatestStart{\mSequence'} = \mNodeDue{i}$, and $\mRefTWFeasible{\mSequence'} = \textsc{true}$. A route $\mRoute{}$ satisfies the time window constraint if for any two subsequences $\mSequence^1$ and $\mSequence^2$, $\mRefEarliestCompletion{\mSequence^1} + \mArcTime{\mSequence^1}{\mSequence^2} \leq \mRefLatestStart{\mSequence^2}$ holds, i.e., if we can reach the first visit of $\mSequence^2$ after serving $\mSequence^1$ exactly at or before the latest start time of $\mSequence^2$.

\textit{Cost evaluation:} Let $\mRefCost{\mSequence}$ be the accumulated cost of sequence $\mSequence$. Then, we define our concatenation operation $\oplus$ as:
\begin{align}
    \mRefCost{\mSequence^1 \oplus \mSequence^2} &= \mRefCost{\mSequence^1} + \mArcCost{\mSequence^1}{\mSequence^2} + \mRefCost{\mSequence^2}
\end{align}

\textit{Efficiency:} Each concatenation operation detailed above can be executed in \bigO{1} time. $\mRefFw{i}$ and $\mRefBw{i}$ can thus be calculated in \bigO{|2\mSetNodes|} by traversing the route using $\mSucc{i}$, starting from the vehicle node $\mVehicleNode{}$, to perform $\mRefFw{\mSucc{i}} = \mRefFw{i} \oplus i$, and $\mPred{i}$, starting from the last visit $\mRoutePosition{\mCardinality{\sigma}}$, to set $\mRefBw{\mPred{i}} = i \oplus \mRefBw{i}$, respectively.
When performing any changes to a route, $\mRefFw{i}$ and $\mRefBw{i}$ need to be updated accordingly. 

For each vehicle $\mVehicleNode_{j} \in \mSetVehicles$, we additionally maintain $\mLast{\mVehicleNode_{j}}$ which holds the information on the last visit of route $\mRoute{j}$. The objectives can be obtained as follows: i) $|\{i\in \mSetNodes: \mVehicleAssignment{i} = \mCardinality{\mSetVehicles}+1\}|$ is the number of unassigned requests, and ii) $\mSum{j \in \mSetVehicles}{} \mRefCost{\mRefFw{\mLast{j}}}$ the total travel cost.

\subsubsection{Iterative local search}
\label{section:methodology:ils}
Algorithm~\ref{alg:decomposition} shows the pseudocode of our \gls{ils}-based metaheuristic. First, we generate an initial solution using a simple but fast parallel insertion construction heuristic (l.~\ref{alg:decomposition:init}). Herein, we randomly assign one request to each vehicle and distribute the remaining requests sequentially to their best position in any route.
\begin{algorithm}[tbh]
\caption{Iterative local search approach}
\label{alg:decomposition}
\SetKwData{TimeLimit}{$\mTimeLimit$}
\SetKwData{Threshold}{$\mThreshold$}
\SetKwData{ThresholdInit}{$\mThresholdInit$}
\SetKwData{ThresholdDec}{$\mThresholdDec$}
\SetKwData{Solution}{$\mSolution$}
\SetKwData{PartSolutions}{$\mPartialSolutions$}
\SetKwData{PartSolution}{$\mSolution_i$}
\SetKwData{CurrentSolution}{$\mSolution'$}
\SetKwData{NewSolution}{$\mSolution''$}
\SetKwData{BestSolution}{$\mSolution^*$}
\SetKwData{AvgSplitSize}{$\mAvgSplitSize$}
\SetKwData{MaxIterations}{$\mMaxIterations{A}$}
\SetKwData{MaxIterationsSubParts}{$\mMaxIterationsSplit$}
\SetKwData{PerturbationSize}{$\mPerturbationSize{A}$}
\SetKwData{Blocks}{$\mBlocks$}
\SetKwData{Graph}{$\mGraphKDSP$}
\SetKwData{KDSP}{$\mKDSP$}
\SetKwFunction{InitialSolution}{parallel\_insertion\_construction}
\SetKwFunction{Decompose}{decompose}
\SetKwFunction{RuinAndRecreate}{ruin\_and\_recreate}
\SetKwFunction{Recombine}{recombine}
\SetKwFunction{AcceptSolution}{accept\_solution\_rnr}
\SetKwFunction{AcceptSolutionILS}{accept\_solution\_ils}
\SetKwFunction{Perturb}{perturb}
\SetKwFunction{ExtractBlocks}{extract\_blocks}
\SetKwFunction{CreateGraph}{create\_graph}
\SetKwFunction{RunKDSP}{run\_kdsp}
\SetKwFunction{BuildFromKDSP}{build\_from\_kdsp}
\KwIn{Time limit \TimeLimit; Average split size \AvgSplitSize; Maximum iterations \MaxIterations; Maximum iterations per split solution \MaxIterationsSubParts; Maximum perturbation size \PerturbationSize}
$\Solution \leftarrow \InitialSolution()$\;\label{alg:decomposition:init}
$\CurrentSolution, \BestSolution \leftarrow \Solution$\;
$i^{\textsc{ils}} \leftarrow 0$\;\label{alg:decomposition:init-total-iterations}
\While{$\TimeLimit \text{ not reached}$}{\label{alg:decomposition:main-loop}
    $i^{\textsc{ls}} \leftarrow 0$\;\label{alg:decomposition:init-iterations}
    $\Threshold \leftarrow \ThresholdInit$\;
    \While{$i^{\textsc{ls}} < \MaxIterations$ \label{alg:recombination:loop_begin}}{
        $\PartSolutions \leftarrow \Decompose(\CurrentSolution, \AvgSplitSize)$\;\label{alg:decomposition:decompose}
        \For{$\PartSolution$ $\textbf{in}$ $\PartSolutions$}{
            $\PartSolution \leftarrow \RuinAndRecreate(\PartSolution, \MaxIterationsSubParts, \Threshold, \ThresholdDec)$\;\label{alg:decomposition:rnr}        
        }
        $\Blocks \leftarrow \ExtractBlocks(\PartSolutions)$\;\label{alg:recombination:extract_blocks}
        $\Graph \leftarrow \CreateGraph(\Blocks)$\;\label{alg:recombination:create_graph}
        $\KDSP \leftarrow \RunKDSP(\Graph)$\;\label{alg:recombination:run_kdsp}
        $\NewSolution \leftarrow \BuildFromKDSP(\KDSP)$\;\label{alg:recombination:build_from_kdsp}
        $\Solution,\CurrentSolution,\BestSolution \leftarrow \AcceptSolution(\NewSolution, \CurrentSolution, \Solution,\BestSolution)$\;\label{alg:decomposition:accept-solution-inner}
        $\Threshold \leftarrow \Threshold - \ThresholdDec \MaxIterationsSubParts$\;\label{alg:decomposition:update-acceptance-criterion}
        $i^{\textsc{ls}} \leftarrow i^{\textsc{ls}} + \MaxIterationsSubParts$\; \label{alg:decomposition:inner_inc}\label{alg:recombination:loop_end}
    }
    $\BestSolution, \CurrentSolution \leftarrow \AcceptSolutionILS(\NewSolution, \CurrentSolution, \BestSolution)$\;\label{alg:decomposition:accept-solution-outer}
    $\CurrentSolution \leftarrow \Perturb(\CurrentSolution, \PerturbationSize)$\;\label{alg:decomposition:perturb}
    $i^{\textsc{ils}} \leftarrow i^{\textsc{ils}} + 1$\; \label{alg:decomposition:inc}
}
\Return $\Solution$
\end{algorithm}

Then, we run the main \gls{ils} loop until we reach a time limit $\mTimeLimit$ (l.~\ref{alg:decomposition:main-loop}). For each \gls{ils} iteration, we perform our local search for $\mMaxIterations{A}$ iterations (l.\ref{alg:recombination:loop_begin}--l.\ref{alg:recombination:loop_end}): herein, we first decompose the problem by inspecting the current solution $\mSolution'$ and randomly partitioning its routes such that each part contains $\mAvgSplitSize$ requests on average. 
Then, we distribute unassigned requests between the partitions. 
To do so, we maintain a history comprising two counters tracking i) how often requests are scheduled directly after another and ii) how often requests are assigned to the same route. For each unassigned request, we sum up the counters for each request in a partition and use this sum as weight in a roulette-wheel selection.

Each resulting partition constitutes a subproblem containing the respective vehicle and requests, for which we use the vehicle's sequence in $\mSolution'$ to obtain an initial solution. We then improve these solutions in parallel by using an \gls{rnr} metaheuristic (l.~\ref{alg:decomposition:rnr}). After a predefined number of $\mMaxIterations{S}$ iterations, we collect the best solution found for each subproblem and recombine these. 
Herein, we first identify and extract the block information from the partial solutions (l.~\ref{alg:recombination:extract_blocks}). For each block, i.e., a sequence of pickup and delivery visits, we obtain the total time spent (service time, waiting time, and travel time), the total distance traveled, as well as the earliest and latest beginning of service time, such that the waiting time is minimal. Note that the block starting time window depends on the partial solutions, resulting in tighter bounds than if computed in isolation. 

Next, we create the dispatching graph $\mGraphKDSP= (\mSetNodesKDSP,\mSetArcsKDSP)$ (l.~\ref{alg:recombination:create_graph}), analog to Section~\ref{section:methodology:matheuristic:dispatching}. After solving the \gls{kdspp} in line~\ref{alg:recombination:run_kdsp}, a set of disjoint shortest paths $\Pi$ has been produced, which is used in line~\ref{alg:recombination:build_from_kdsp} to generate routes by concatenating the sequences of the corresponding routes until reaching the sink.

During the search, we may allow a decrease in quality when accepting incumbent solutions to balance exploration and intensification (l.~\ref{alg:decomposition:accept-solution-inner}). 
Here, we adopted the record-to-record aspiration criterion, which performs well for the classical \gls{pdptw} instances \citep[see, e.g.,][]{SantiniRopkeEtAl2018}, as follows. New best solutions are always accepted; otherwise, a solution is only accepted if the gap between it and the current best solution is smaller than a threshold $\mThreshold$, which decreases linearly by $\mThresholdDec$ with each iteration. We use $\mThresholdInit=0.333$ and $\mThresholdDec=\mThresholdInit / \mMaxIterations{}$ as parameters for the initial temperature and decrement step-value. The \gls{rnr} procedures run in parallel for each subproblem and start with the current criterion state. We update the temperature after recombining the solutions by decreasing the temperature $\mMaxIterations{S}$ times (l.~\ref{alg:decomposition:update-acceptance-criterion}).

After the local search procedure finishes, we perform a second check for accepting the solution found (l.~\ref{alg:decomposition:accept-solution-outer}). We accept any new best solution found so far. Otherwise, we either use the current solution found by the cost-minimization component or revert to the best solution. This choice is probability-driven and controlled by two counters, the total number of \gls{ils} iterations and the number of \gls{ils} iterations $i^{\textsc{ils}}$ without any improvements to the best solution $i^{\textsc{last}}$, dividing the latter by the former, $i^{\textsc{last}}/i^{\textsc{ils}}$.

The perturbation procedure attempts to escape local optima (l.~\ref{alg:decomposition:perturb}) by performing either a random relocation or an exchange operation with a predefined probability $\mProbRelocateExchange$ for a total number of $Z_A$ times. This random relocation moves a request to a feasible position in another route, which is again chosen randomly: the exchange operations randomly select two requests in distinct routes and attempt to insert them in the other routes at a random feasible location.

\begin{figure}[t]
    \centering
    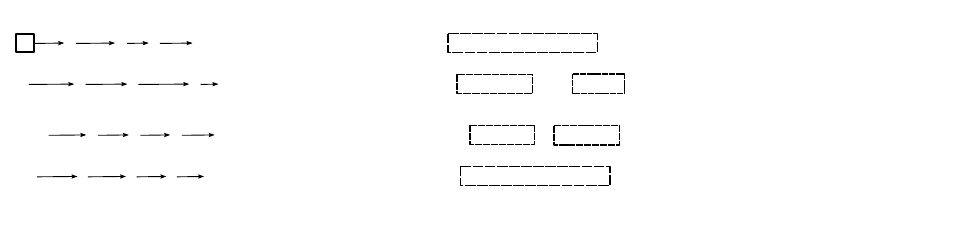
    \caption{Example of a run of the recombination procedure.}
    \begin{subfigure}{0.28\textwidth}
    \subcaption{Partial solutions $\mSolution^1$ and $\mSolution^2$.}
    \label{fig:block_recombination:parts}
    \end{subfigure}
    \begin{subfigure}[c]{0.05\textwidth}\end{subfigure}
    \begin{subfigure}[c]{0.435\textwidth}
    \subcaption{Auxiliary graph $\mGraphKDSP$ with the k-dSPP arcs in solid black.}
    \label{fig:block_recombination:auxgraph}
    \end{subfigure}
    \begin{subfigure}[c]{0.05\textwidth}\end{subfigure}
    \begin{subfigure}[c]{0.25\textwidth}
    \subcaption{Resulting recombined complete solution.}
    \label{fig:block_recombination:result}
    \end{subfigure}

    \label{fig:block_recombination}
\end{figure}

\subsubsection{Ruin-and-recreate procedure}\label{section:methodology:rnr}

The core component of our algorithm follows the \gls{rnr} approach \citep{SchrimpfSchneiderEtAl2000}. In each iteration, parts of the solution are ruined (l.~\ref{alg:rnr:ruin}) by removing requests and then recreated (l.~\ref{alg:rnr:recreate}) by reinserting the unassigned requests again in a cost-improving manner. We follow the implementation of the slack-inducing string removal approach from \cite{ChristiaensVandenBerghe2020}. Herein, the ruin operator is an adjacent string removal operator, which selects an assigned request randomly and then iteratively removes sequences of visits from routes containing related requests to introduce slack for potentially better re-insertions. 
The recreate operator is a variant of a simple greedy best insertion approach. Unassigned requests are sorted and then sequentially inserted to their best position in any route. A so-called blink mechanic adds diversity to this approach, wherein insertion positions may be skipped with a slight chance. We modified this approach to test the insertion of at most 40 requests per operator call to reduce the computational runtime. 

Finally, the \gls{rnr} procedure stops after a predefined number of $\mMaxIterations{A}$ iterations has been performed. In the following, we detail the operators and intensification methods applied.

\begin{algorithm}[hbt]
\caption{Ruin and Recreate}
\label{alg:rnr}
\SetKwData{Solution}{$\mSolution$}\SetKwData{CurrentSolution}{$\mSolution'$}\SetKwData{NewSolution}{$\mSolution''$}
\SetKwData{MaxIterations}{$M_A$}\SetKwData{Threshold}{$\mThreshold$}
\SetKwData{ThresholdInit}{$\mThresholdInit$}\SetKwData{ThresholdDec}{$\mThresholdDec$}
\SetKwFunction{Ruin}{ruin}
\SetKwFunction{Recreate}{recreate}
\SetKwFunction{LS}{local\_search}
\SetKwFunction{AcceptSolution}{accept\_solution}
\SetKwFunction{Rand}{random}
\KwIn{Feasible solution \Solution; Maximum iterations \MaxIterations; Inital acceptance threshold \Threshold; Theshold decrement value \ThresholdDec}
$\CurrentSolution \leftarrow \Solution$\;
$\Threshold \leftarrow \ThresholdInit$\;
\For{$i\leftarrow 1$ \KwTo $\MaxIterations$}{
$\NewSolution \leftarrow \Ruin(\CurrentSolution)$\;\label{alg:rnr:ruin}
$\NewSolution \leftarrow \Recreate(\NewSolution)$\;\label{alg:rnr:recreate}
$\Solution, \CurrentSolution \leftarrow \AcceptSolution(\Solution, \CurrentSolution, \NewSolution, \Threshold)$\;\label{alg:rnr:accept-solution}
$\Threshold \leftarrow \Threshold - \ThresholdDec$\;
}
\Return $\Solution$
\end{algorithm}

\paragraph{Ruin: Adjacent string removal.}

The ruin operator starts by selecting a request at random, which we refer to as $\mRequest{}^*$. Let $\overline{\mSetRequests}(\mRequest{}^*)$ be a sorted set of related requests of $\mRequest{}^*$, including $\mRequest{}^*$, and let $\mVehicleAssignment{\mRequest{i}}$ be the current vehicle assignment of the respective pickup node of $\mRequest{i}$. 
For every $\mRequest{i} \in \overline{\mSetRequests}(\mRequest{}^*)$, the algorithm selects a consecutive sequence of $k_s$ visits from route $\mVehicleAssignment{\mRequest{i}}$ and either a) and removes the corresponding requests (called split), or b) selects a substring $\bar{{\mRoute{}}} \subset {\mRoute{}}$, which is ignored, and remove the requests of the remaining consecutive string (called split-string). This process continues with the following requests, skipping routes already modified and terminating if a certain threshold of routes $k_s$ has been considered. 

Both the number of strings to remove $k_s$ and the length of the string $l_{\mRoute{}}$ are randomly chosen from a uniform distribution. We bound them using two parameters: i) an average number of nodes to remove $\bar{c}$, and ii) a maximum number of removed strings $L^{max}$, using the following equations
\begin{align}
    k_s &=  \left\lfloor \mathcal{U}(1,\frac{4\bar{c}}{1+l_s^{max}} - 1) \right\rfloor \\
    l_{\mRoute{}} &=  \left\lfloor \mathcal{U}(1,\text{min}\{ \left|{\mRoute{}}\right|, l_s^{max}\}) \right\rfloor
\end{align}
where, $l_s^{max} = \text{min}\{ L^{max}, \left|\overline{{\mRoute{}} \in \mSetVehicles}\right|\}$ holds the maximum string cardinality, which is the minimum of the parameter $L^{max}$ and the average cardinality of the current routes. We use $\bar{c} = 15$ and $L^{max} = 10$ as identified in our parameter tuning (see Appendix~\ref{appendix:parameter-tuning-details}).

We randomly perform the split with probability $\mProbRuinMode=0.75$, or split-string sub-procedure otherwise, aiming to remove $l_{\mRoute{}}$ nodes. The split procedure selects a consecutive string of nodes, which are removed. In split-string, we first calculate the length of the substring $m$, where $l_{\mRoute{}} + m \leq \mCardinality{{\mRoute{}}}$. We start with $m=1$ and increment $m$ by one with probability $\mRuinBeta=0.10$ until it has either not been incremented or reached a maximum of $\mCardinality{{\mRoute{}}} - l_{\mRoute{}}$. After removal, we check whether requests were only partially removed (either pickup or delivery node) and fully remove them from the solution.

\paragraph{Recreate: Greedy insertion with blinks.}

To recreate our solution, we use a greedy insertion heuristic. The procedure sequentially tries to insert pickup-and-delivery-pairs, i.e., requests, to their best possible positions, ensuring they satisfy the precedence constraint. 
Similar to \cite{ChristiaensVandenBerghe2020}, we ignore some insertion positions while searching for the best solution with a low probability. These so-called blinks add randomness to the insertion procedure similar to other noise-based methods from the literature, e.g., \cite{RopkePisinger2006}.

The order in which requests are inserted is determined anew every time the procedure is called using a weighted roulette-wheel selection. We use six different criteria (random, far, close, tw-length, tw-start, tw-end) and their corresponding weights (6, 2, 1, 4, 2, 2) to order the unassigned requests based on their related pickup nodes i) at random, ii) descending or iii) ascending by their distance to the closest vehicle starting location, iv) ascending by the time window length, v) ascending by the earliest start time, or vi) descending by the latest start time.

The routes are searched sequentially. First, routes already serving requests are considered. If no insertion position is found -- either because it cannot be inserted feasibly or was randomly disregarded by the blinking mechanic -- all currently empty routes are considered for insertion, and one with the minimum additional cost is selected.

When traversing a route to look for a feasible insertion point, we consider insertion positions by first considering the pickup node position and then testing all viable subsequent positions for inserting the delivery node. By doing so in a lexicographical order, we can efficiently use the preprocessed evaluation data $\mRefFw{i}$ and $\mRefBw{j}$, as follows.

When testing an insertion of request $r=(i, j)$ after node $i'$ and before node $j' = \mSucc{i'}$, we can use the preprocessed data defined in Section~\ref{section:methodology:eval}. Figure~\ref{fig:insertion:preprocessed} shows the sequences after insertion if done in a lexicographical manner. The corresponding concatenation operations are shown on the right side using the preprocessed route data. We can observe that the first part of the sequence up to but excluding the delivery node $\mDeliveryNode$ is part of the following sequence. By reusing this data, we can test all pickup-and-delivery positions efficiently in $\bigO{\mCardinality{\mRoute{}}^2}$.

\begin{figure}[h]
    \centering
    \begin{tabular}{ll}
    \toprule
    Sequence after insertion & Evaluation using preprocessed data\\
    \midrule
    $\mSequence' = \{\ldots,i,\mPickupNode,\mDeliveryNode,j,j',j''\ldots\}$ & $\mRefFw{i'} \oplus \mPickupNode \oplus \mDeliveryNode \oplus \mRefBw{j'}$\\
    $\mSequence'' = \{\ldots,i,\mPickupNode,j,\mDeliveryNode,j',j''\ldots\}$ & $\mRefFw{i'} \oplus \mPickupNode \oplus j \oplus \mDeliveryNode \oplus \mRefBw{j'}$\\
    $\mSequence''' = \{\ldots,i,\mPickupNode,j,j',\mDeliveryNode,j''\ldots\}$ & $\mRefFw{i'} \oplus \mPickupNode \oplus j \oplus j' \oplus \mDeliveryNode \oplus \mRefBw{j''}$\\
    \vdots\\
    \bottomrule
    \end{tabular}
    \caption{Sequences considered when inserting request $r=(\mPickupNode,\mDeliveryNode)$ between nodes $i,j$ and the concatenation operations to evaluate them.}
    \label{fig:insertion:preprocessed}
\end{figure}

\paragraph{Accepting solutions.}

During the search, we may allow a decrease in quality when accepting incumbent solutions to balance exploration and intensification (l.~\ref{alg:rnr:accept-solution}). Herein, we follow the same linear record-to-record approach described in Section~\ref{section:methodology:ils}. 
When encountering a better solution, we always accept it and continue the search.

\paragraph{Intensification procedure.}

To support intensification, we perform an intra-route local search approach based on \cite{BalasSimonetti2001} to 
find efficient itineraries whenever we reach a new best solution. Herein, a restricted neighborhood with visits up to a predefined number of positions, $k$, from their current location can be exhaustively searched efficiently. 

Let $k$ be a positive integer and $\mRoute{}$ be a feasible initial route, where $\mRoutePosition{u}$ denotes the position of node $u$ in route $\mRoute{}$. The BS($k$) neighborhood is defined as all permutations of visits $\pi$ in $\mRoute{}$, satisfying
\begin{align}
    \mPermutationPosition{u} &< \mPermutationPosition{v} && \forall u,v \in \mSetNodes: \mRoutePosition{u} + k \leq \mRoutePosition{v} \label{eq:bs:k}
\end{align}
where $\mPermutationPosition{u}$ denotes the position of node $u$ in the permutated sequence. 
Equation \eqref{eq:bs:k} enforces that permutations respect the maximum number of positions moved by any node, $k$, from their original position in $\mRoute{}$. 
Similar to \cite{BalasSimonetti2001} the search for the best permutation $\mPermutation$ satisfying (\ref{eq:bs:k}) reduces to a shortest path problem on an auxiliary acyclic graph $G^*=(V^*, E^*)$. This graph is arranged into layers so that only successive layers are connected by arcs. There are $2\left|\mSetRequests\right|+2$ layers required at maximum. The first layer and last correspond to the departure of a vehicle from its starting location, and a dummy node represents the end of the route. The other layers correspond to the $2\left|\mSetRequests\right|$ visits in $\mPermutation$. Each layer contains multiple nodes representing distinct states of anticipated and delayed visits to locations in relation to $\mRoute{}$. More specifically, each node in $V^*$ represents a quadruplet $(i,i+j,S^-(\mRoute{},i),S^+(\mRoute{},i))$, where $i$ is the current position in $\mPermutation$; $j$ is the offset from the current location, which visit allocated to position $i$ in $\mPermutation$, such that $\mPermutationPosition{i} = \mRoutePosition{i+j}$; $S^-(\mRoute{},i)$ is the set of visits that are allocated at position $i$ or later in $\mRoute{}$ and allocated at position $i$ or before in $\mPermutation$; and, finally, $S^+(\mRoute{},i)$ is the set of visits allocated before position $i$ in $\mRoute{}$ that will be allocated at position $i$ or after in $\mPermutation$.

Observe that $S^-(\mRoute{},i)$ and $S^+(\mRoute{},i)$ provide enough information to restrict the search to solutions satisfying Constraints~\eqref{eq:bs:k}, by restricting the search to nodes such that all visits in $S^+(\mRoute{},i)$ are located in no more than $k$ positions earlier than $v$ in $\mRoute{}$. The precedence constraint can be trivially checked by maintaining the information of allocated pickup nodes in each node and prohibiting traversing invalid arcs. 

\begin{figure}[bt]
    \centering
    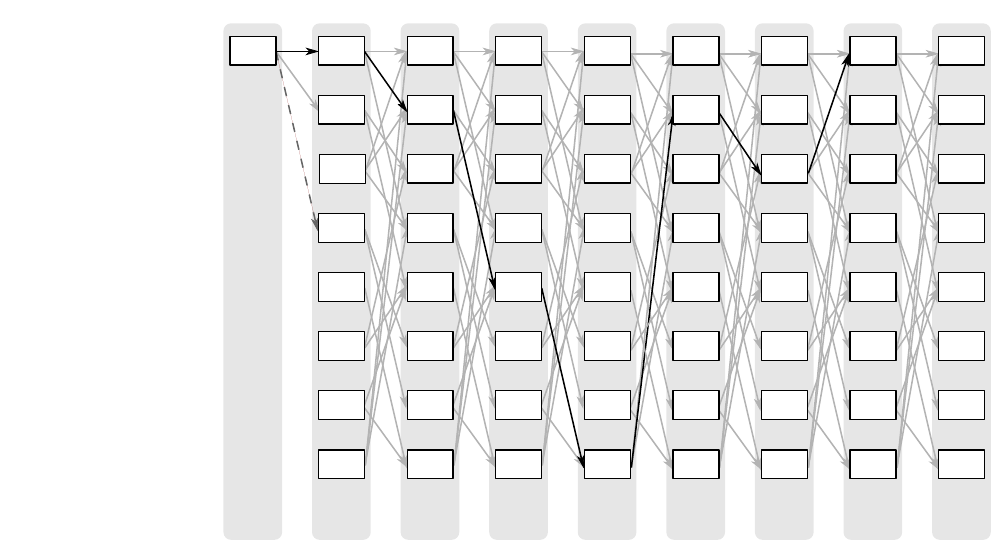
    \caption{Graph $G^*$ with $k=3$ and $\mRoute{}=\{\mVehicleNode{}, \mPickupNode_{1}, \mPickupNode_{2}, \mDeliveryNode_{1}, \mPickupNode_{3}, \mDeliveryNode_{2},\mDeliveryNode_{3}\}$ (Nodes contain the corresponding pickup or delivery node for visualization). Arcs in black show the path of permutation $\pi = \{\mVehicleNode{}, \mPickupNode_{1}, \mDeliveryNode_{1}, \mPickupNode_{3}, \mPickupNode_{2}, \mDeliveryNode_{3}, \mDeliveryNode_{2}\}$. Note that the dashed arc from Layer 1, Node 1 to Layer 2, Node 4 is invalid as the resulting sequence would violate the precedence constraint.}
    \label{fig:methodology:bs}
\end{figure}

Figure \ref{fig:methodology:bs} shows $G^*$ on a simple example with three requests, starting with route $\mRoute{}=\{\mVehicleNode{}, \mPickupNode_{1}, \mPickupNode_{2}, \mDeliveryNode_{1}, \mPickupNode_{3}, \mDeliveryNode_{2},\mDeliveryNode_{3}\}$ and $k=3$. 
Without considering capacity or time window constraints, finding the shortest path in this graph can be done in $\mathcal{O}(k^22^{k-2}\mCardinality{\mRoute{}})$ time. 
Note that graph $G^*$ need only to be built once -- requiring $\mathcal{O}(k2^{k-2}n)$ space -- and can be reused. 
Furthermore, by fixing $k$, we obtain a polynomial-time algorithm to locate the best permutation.

Time window constraints cannot be included in the graph structure and need to be checked during the search, rendering the problem of finding the shortest path with resource constraints. This problem is NP-hard \citep{IrnichDesaulniers2005} and requires a labeling algorithm that no longer runs in polynomial time. Each vertex in the graph now may contain multiple labels, each representing a non-dominated path through the graph $G^*$. We reuse the evaluation mechanics from Section~\ref{section:methodology:eval} to define a label $L(v^*, \mSequence')=(\mRefCost{\mSequence'}, \mRefCapacitySum{\mSequence'}, \mRefCapacityMax{\mSequence'}, \mRefTraveltime{\mSequence'}, \mRefEarliestCompletion{\mSequence'}, \mRefLatestStart{\mSequence'})$, where $v^*$ is the corresponding vertex in $G^*$ and $\mSequence$ the sequence of visits reaching and including the node added at $v^*$. Note that each label at any vertex $v^*$ visited the same nodes in $\mSetNodes$ but in a different order.

When extending a label $L$ to a compatible successor, we first check whether the new sequence is feasible regarding the capacity and time window constraint. If so, we create a new label for the extended sequence, add it to the corresponding label set, and perform a dominance check. Given two labels $L^1(v^*,\mSequence^1)$ and $L^2(v^*,\mSequence^2)$, we say $L^1$ dominates $L^2$ if all the following hold:
\begin{align}
    \mRefCost{\mSequence^1} &\leq \mRefCost{\mSequence^2} \label{eq:dom:cost} \\
    \mRefEarliestCompletion{\mSequence^1} &\leq \mRefEarliestCompletion{\mSequence^2} \label{eq:dom:completion}
\end{align}

While Equation \eqref{eq:dom:cost} ensures that we keep the lowest cost label, \ref{eq:dom:completion} ensures that a label where we arrive earlier stays as well, as further extensions may invalidate labels with later arrivals due to violations in the time window constraint.
As the number of labels per vertex can grow exponentially due to the potential trade-off in cost and time, we limit the number to only consider the four least-cost labels. Note that this adaption may lead to searches terminating without finding any feasible solution.

\section{Experimental Design}\label{section:experimental-design}

In this section, we first introduce our case-study data and instances. We then generate mid-sized benchmark instances for our algorithmic component analysis.

We implemented the algorithms presented in Section~\ref{section:methodology} in the Rust programming language, compiled using rustc version 1.70 in release mode and with link-time optimization. The experiments were conducted on a cluster setup consisting of Intel i9-9900 processors with 16 cores running at 3.10 GHz and 64GB RAM, using Ubuntu 20.04 LTS operating system. We limit our runs to four cores and 8GB of RAM to allow conducting our experiments in parallel. Furthermore, we run experiments with five different seeds to reduce the probability of outliers.

The source code, results, and documentation on building and running the algorithm are available online at \url{https://github.com/tumBAIS/LargeScalePDPTW}.

\subsection{Case-study instances}\label{subsection:case-study}

\begin{figure}
    \centering
     \begin{subfigure}[b]{0.36\textwidth}
         \centering
            \includegraphics[height=6cm]{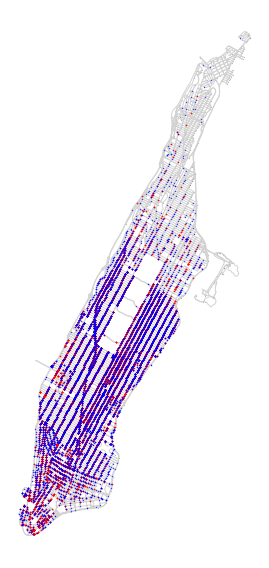}

            \caption{New York City street network from OSM with pickups and dropoff points marked in red and blue, respectively.}
            \label{fig:study:streetnetwork}
     \end{subfigure}
     \hfill
     \begin{subfigure}[b]{0.56\textwidth}
         \centering
            {
                \small
                \input{figures/requests_arriving_rolling_10min_all_days}
            }
            \caption{Arriving trips per minute during the sample week (10-minute rolling average). Selected time frames are highlighted in green.}
            \label{fig:study:time-frames}
     \end{subfigure}
     \caption{New York City trip and network data.}
\end{figure}

One of the largest and most regularly used data sources in the literature is the data from the NYC Taxi and Limousine Commission \citep{NYC_TLC}, which provides historical taxi trip data. We use a preprocessed version of these trip data by \cite{JungelParmentierEtAl2023}, available online at \url{https://github.com/tumBAIS/ML-CO-pipeline-AMoD-control}. This data includes origin and destination location information as well as timestamps. We use \gls{osm} network data to map the trip data and calculate travel time and distances, assuming an average travel speed of 20 km/h. Pickup and dropoff locations are matched to the nearest node of the street network (see Figure~\ref{fig:study:streetnetwork}). The distance and travel time matrix is preprocessed once for the whole network using an all-pair shortest path algorithm. 

Figure~\ref{fig:study:time-frames} shows the 10-minute rolling average of trips arriving throughout the day in our sample week, spanning January 5\textsuperscript{th} to January 11\textsuperscript{th} in 2015. We observe similar patterns for weekdays (Mon-Fri) and diverging patterns for the weekend days (Sat, Son), e.g., at 2 a.m. For this reason, we generate our large-scale instances only for weekdays and extract trip data for six distinct time frames for one hour each: night (2-3 a.m.), early morning (6-7 a.m.), late morning (10-11 a.m.), afternoon (2-3 p.m.), late afternoon (6-7 p.m.), late evening (10-11 p.m.).

Requests are generated from the trip data of a given day and time frame. We use the dropoff time of the trip data as the fixed point and calculate the corresponding (earliest) pickup time based on the direct travel time, i.e., $e_{\mPickupNode} = e_{\mDeliveryNode} - \mDirectArcTime{}$. The latest begin-of-service times $l_\mPickupNode, l_\mDeliveryNode$ are implicitly defined by the buffer $\mBuffer$. 

Vehicles start at random trip destinations, sampled from the full data set. Capacities of vehicles are set to three requests when pooling and one for taxi operations. We consider a base scenario with a fleet of 1000 vehicles. 
Table~\ref{tab:nyc_instances} shows the number of requests in the generated instances.

\begin{table}[htb]
    \centering
    \begin{tabular}{cccccccc}
    \toprule
        Date & Day & 2-3h & 6-7h & 10-11h & 14-15h & 18-19h & 22-23h \\
     \midrule
     2015-01-05 & Mon. & 1.385 & 6.742 & 11.457 & 13.094 & 18.736 & 10.070 \\
     2015-01-06 & Tue. & 1.035 & 7.085 & 11.922 & 13.441 & 19.490 & 12.150 \\
     2015-01-05 & Wed. & 1.669 & 7.296 & 14.086 & 15.256 & 21.375 & 15.633 \\
     2015-01-08 & Thu. & 1.815 & 7.705 & 15.126 & 15.212 & 21.170 & 17.115 \\
     2015-01-09 & Fri. & 2.869 & 6.878 & 13.552 & 14.820 & 21.222 & 18.691 \\
    \bottomrule
    \end{tabular}
    \caption{Number of requests for the instances generated from the NYC data set.}
    \label{tab:nyc_instances}
\end{table}

\begin{figure}[tbh]
    \centering
    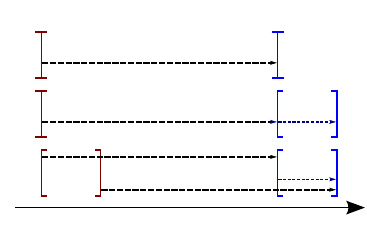
            \caption{Time window settings considered in the case study (A: fixed pickup and dropoff; B: fixed pickup, variable dropoff; C: variable pickup and dropoff).}
            \label{fig:study:timewindow}
\end{figure}

Our experiments are designed to capture the marginal benefits of transitioning from classic taxi operations to flexible ride-pooling offerings; see Figure~\ref{fig:study:timewindow}.
In the first setting (A), we do not allow any deviation from the desired pickup and dropoff time, emulating the current operations of taxi fleets. Nevertheless, we consider ride-pooling in this setting to identify (if any) on-route pooling opportunities. In the second setting (B), we consider deviating from the desired dropoff time by $\mBuffer$ minutes, $\mBuffer=\{1,2,3,4,6,8,10,12,15,20\}$. This does not impact taxi operations as the pickup time is still fixed. However, ride-pooling operations may benefit from this increased flexibility, allowing more requests to be served concurrently. In setting (C), which is the original problem setting (see Section~\ref{section:problem_formulation}), we extend the time windows for pickup and dropoff, allowing for even more flexibility on the operator side. For this setting, we conduct our experiments with taxi operations and ride-pooling, as the former may also benefit from the flexibility.

Note that solutions found for a buffer setting $\mBuffer_a$ are feasible for $\mBuffer_b$ where $\mBuffer_a < \mBuffer_b$. We exploit this property in our managerial experiments as follows. We start by solving the most restrictive case of zero buffer ($\mBuffer=0)$ with a runtime of 8 hours for each day and time frame. All subsequent settings with $\mBuffer\geq0$s are run for 4 hours but warm-started using the solution of the previous buffer setting. 
With this approach, we avoid issues arising from the heuristic nature of the algorithm. 
Note further, that for taxi operations with no flexibility, we only perform a single run per instance, as the \gls{kdspp} always finds the optimal solution in this setting.

\subsection{Mid-sized benchmark set}

We create two additional sets comprising mid-sized instances following the instance generation protocol detailed in Section~\ref{subsection:case-study} -- one for our computational study comparing the sequential and integrated approach and one for algorithmic tuning and analysis.
In both cases, we consider a subset of requests for the same time horizon of one hour. To ensure the temporal properties of the large-scale instances are retained, we do not select a subset of requests at random but equally distanced.

Our first set, $\textsc{Mid}^{25}$, comprises the period 18-19h only, where we select 25\% of the original requests. These instances comprise around 5000 requests. Here, the aim is to compare the sequential, integrated, and combined approaches on high-demand instances. We consider buffer values from $\delta=1,\ldots,6$ minutes to analyze their performances on varying degrees of freedom. Furthermore, we compare three fleet size configurations: a) a low availability setting with 200 vehicles, b) availability of 500 vehicles, where we have roughly one vehicle per 10 requests, and c) high availability of 800 vehicles.

The second set is created for analysis and tuning of the metaheuristic approach. For set $\textsc{Mid}^{10}$, we select 10\% of the original requests and omit the early morning periods of 2–3 a.m. and 6–7 a.m. due to their small size. We consider all buffer settings as in the large instances and generate three different fleet size configurations with 100, 200, and 300 vehicles, respectively. From this setting, we randomly select one instance per day, resulting in 18 instances in total.
\section{Results}\label{section:results}
In the following, we discuss our numerical results. First, we focus on a pure algorithmic perspective and compare the performance of our proposed algorithms: the sequential matheuristic ($\mSolverSequentialMathHeur$), the iterated local search that takes integrated decisions ($\mSolverDBILS$), and our $\mSolverDBILS$ when using the solution of $\mSolverSequentialMathHeur$ as a warmstart ($\mSolverCombined$). Second, we compare our developed algorithm against a known benchmark data set from the literature that provides instances for the \gls{pdptw} in a medium to large-scale ride-hailing context. Finally, we use our best-performing algorithm to derive managerial insights for the presented large-scale case study.

\subsection{Computational Analyses}

\begin{figure}[!htb]
    \centering
     \begin{subfigure}[b]{0.48\textwidth}
        \centering
\begin{tikzpicture}

\definecolor{darkslategray38}{RGB}{38,38,38}
\definecolor{darkslategray61}{RGB}{61,61,61}
\definecolor{lightgray204}{RGB}{204,204,204}
\definecolor{peru22412844}{RGB}{224,128,44}
\definecolor{seagreen5814558}{RGB}{58,145,58}
\definecolor{steelblue49115161}{RGB}{49,115,161}

\begin{axis}[
width=0.8\textwidth,
height=8cm, 
axis line style={lightgray204},
legend cell align={left},
legend style={fill opacity=0.8, draw opacity=1, text opacity=1, draw=lightgray204, nodes={scale=0.7, transform shape}},
tick align=outside,
x grid style={lightgray204},
xlabel=\textcolor{darkslategray38}{Fleet},
xmajorticks=true,
xmin=-0.5, xmax=2.5,
xtick style={color=darkslategray38},
xtick={0,1,2},
xticklabels={200,500,800},
y grid style={lightgray204},
ylabel=\textcolor{darkslategray38}{Unassigned requests $\mUnassigned$},
ymajorgrids,
ymajorticks=true,
ymin=-154.1, ymax=3236.1,
ytick style={color=darkslategray38}
]
\path [draw=darkslategray61, fill=steelblue49115161, semithick]
(axis cs:-0.397333333333333,2384)
--(axis cs:-0.136,2384)
--(axis cs:-0.136,2729.25)
--(axis cs:-0.397333333333333,2729.25)
--(axis cs:-0.397333333333333,2384)
--cycle;
\path [draw=darkslategray61, fill=peru22412844, semithick]
(axis cs:-0.130666666666667,2039)
--(axis cs:0.130666666666667,2039)
--(axis cs:0.130666666666667,2453.75)
--(axis cs:-0.130666666666667,2453.75)
--(axis cs:-0.130666666666667,2039)
--cycle;
\path [draw=darkslategray61, fill=seagreen5814558, semithick]
(axis cs:0.136,2000)
--(axis cs:0.397333333333333,2000)
--(axis cs:0.397333333333333,2435.75)
--(axis cs:0.136,2435.75)
--(axis cs:0.136,2000)
--cycle;
\path [draw=darkslategray61, fill=steelblue49115161, semithick]
(axis cs:0.602666666666667,270.5)
--(axis cs:0.864,270.5)
--(axis cs:0.864,606.5)
--(axis cs:0.602666666666667,606.5)
--(axis cs:0.602666666666667,270.5)
--cycle;
\path [draw=darkslategray61, fill=peru22412844, semithick]
(axis cs:0.869333333333333,33.5)
--(axis cs:1.13066666666667,33.5)
--(axis cs:1.13066666666667,420.25)
--(axis cs:0.869333333333333,420.25)
--(axis cs:0.869333333333333,33.5)
--cycle;
\path [draw=darkslategray61, fill=seagreen5814558, semithick]
(axis cs:1.136,29)
--(axis cs:1.39733333333333,29)
--(axis cs:1.39733333333333,369)
--(axis cs:1.136,369)
--(axis cs:1.136,29)
--cycle;
\path [draw=darkslategray61, fill=steelblue49115161, semithick]
(axis cs:1.60266666666667,0)
--(axis cs:1.864,0)
--(axis cs:1.864,0)
--(axis cs:1.60266666666667,0)
--(axis cs:1.60266666666667,0)
--cycle;
\path [draw=darkslategray61, fill=peru22412844, semithick]
(axis cs:1.86933333333333,0)
--(axis cs:2.13066666666667,0)
--(axis cs:2.13066666666667,1)
--(axis cs:1.86933333333333,1)
--(axis cs:1.86933333333333,0)
--cycle;
\path [draw=darkslategray61, fill=seagreen5814558, semithick]
(axis cs:2.136,0)
--(axis cs:2.39733333333333,0)
--(axis cs:2.39733333333333,0)
--(axis cs:2.136,0)
--(axis cs:2.136,0)
--cycle;
\draw[draw=darkslategray61,fill=steelblue49115161,line width=0.3pt] (axis cs:0,0) rectangle (axis cs:0,0);
\addlegendimage{ybar,ybar legend,draw=darkslategray61,fill=steelblue49115161,line width=0.3pt}
\addlegendentry{\mSolverSequentialMathHeur}

\draw[draw=darkslategray61,fill=peru22412844,line width=0.3pt] (axis cs:0,0) rectangle (axis cs:0,0);
\addlegendimage{ybar,ybar legend,draw=darkslategray61,fill=peru22412844,line width=0.3pt}
\addlegendentry{\mSolverDBILS}

\draw[draw=darkslategray61,fill=seagreen5814558,line width=0.3pt] (axis cs:0,0) rectangle (axis cs:0,0);
\addlegendimage{ybar,ybar legend,draw=darkslategray61,fill=seagreen5814558,line width=0.3pt}
\addlegendentry{\mSolverCombined}

\addplot [semithick, darkslategray61, forget plot]
table {%
	-0.266666666666667 2384
	-0.266666666666667 2102
};
\addplot [semithick, darkslategray61, forget plot]
table {%
	-0.266666666666667 2729.25
	-0.266666666666667 3212
};
\addplot [semithick, darkslategray61, forget plot]
table {%
	-0.332 2102
	-0.201333333333333 2102
};
\addplot [semithick, darkslategray61, forget plot]
table {%
	-0.332 3212
	-0.201333333333333 3212
};
\addplot [semithick, darkslategray61, forget plot]
table {%
	0 2039
	0 1678
};
\addplot [semithick, darkslategray61, forget plot]
table {%
	0 2453.75
	0 2992
};
\addplot [semithick, darkslategray61, forget plot]
table {%
	-0.0653333333333333 1678
	0.0653333333333333 1678
};
\addplot [semithick, darkslategray61, forget plot]
table {%
	-0.0653333333333333 2992
	0.0653333333333333 2992
};
\addplot [semithick, darkslategray61, forget plot]
table {%
	0.266666666666667 2000
	0.266666666666667 1646
};
\addplot [semithick, darkslategray61, forget plot]
table {%
	0.266666666666667 2435.75
	0.266666666666667 2976
};
\addplot [semithick, darkslategray61, forget plot]
table {%
	0.201333333333333 1646
	0.332 1646
};
\addplot [semithick, darkslategray61, forget plot]
table {%
	0.201333333333333 2976
	0.332 2976
};
\addplot [semithick, darkslategray61, forget plot]
table {%
	0.733333333333333 270.5
	0.733333333333333 143
};
\addplot [semithick, darkslategray61, forget plot]
table {%
	0.733333333333333 606.5
	0.733333333333333 1090
};
\addplot [semithick, darkslategray61, forget plot]
table {%
	0.668 143
	0.798666666666667 143
};
\addplot [semithick, darkslategray61, forget plot]
table {%
	0.668 1090
	0.798666666666667 1090
};
\addplot [black, mark=diamond*, mark size=2.5, mark options={solid,fill=darkslategray61}, only marks, forget plot]
table {%
	0.733333333333333 1164
};
\addplot [semithick, darkslategray61, forget plot]
table {%
	1 33.5
	1 0
};
\addplot [semithick, darkslategray61, forget plot]
table {%
	1 420.25
	1 996
};
\addplot [semithick, darkslategray61, forget plot]
table {%
	0.934666666666667 0
	1.06533333333333 0
};
\addplot [semithick, darkslategray61, forget plot]
table {%
	0.934666666666667 996
	1.06533333333333 996
};
\addplot [semithick, darkslategray61, forget plot]
table {%
	1.26666666666667 29
	1.26666666666667 0
};
\addplot [semithick, darkslategray61, forget plot]
table {%
	1.26666666666667 369
	1.26666666666667 857
};
\addplot [semithick, darkslategray61, forget plot]
table {%
	1.20133333333333 0
	1.332 0
};
\addplot [semithick, darkslategray61, forget plot]
table {%
	1.20133333333333 857
	1.332 857
};
\addplot [black, mark=diamond*, mark size=2.5, mark options={solid,fill=darkslategray61}, only marks, forget plot]
table {%
	1.26666666666667 931
};
\addplot [semithick, darkslategray61, forget plot]
table {%
	1.73333333333333 0
	1.73333333333333 0
};
\addplot [semithick, darkslategray61, forget plot]
table {%
	1.73333333333333 0
	1.73333333333333 0
};
\addplot [semithick, darkslategray61, forget plot]
table {%
	1.668 0
	1.79866666666667 0
};
\addplot [semithick, darkslategray61, forget plot]
table {%
	1.668 0
	1.79866666666667 0
};
\addplot [black, mark=diamond*, mark size=2.5, mark options={solid,fill=darkslategray61}, only marks, forget plot]
table {%
	1.73333333333333 44
	1.73333333333333 55
	1.73333333333333 93
};
\addplot [semithick, darkslategray61, forget plot]
table {%
	2 0
	2 0
};
\addplot [semithick, darkslategray61, forget plot]
table {%
	2 1
	2 1
};
\addplot [semithick, darkslategray61, forget plot]
table {%
	1.93466666666667 0
	2.06533333333333 0
};
\addplot [semithick, darkslategray61, forget plot]
table {%
	1.93466666666667 1
	2.06533333333333 1
};
\addplot [black, mark=diamond*, mark size=2.5, mark options={solid,fill=darkslategray61}, only marks, forget plot]
table {%
	2 34
	2 74
	2 70
};
\addplot [semithick, darkslategray61, forget plot]
table {%
	2.26666666666667 0
	2.26666666666667 0
};
\addplot [semithick, darkslategray61, forget plot]
table {%
	2.26666666666667 0
	2.26666666666667 0
};
\addplot [semithick, darkslategray61, forget plot]
table {%
	2.20133333333333 0
	2.332 0
};
\addplot [semithick, darkslategray61, forget plot]
table {%
	2.20133333333333 0
	2.332 0
};
\addplot [black, mark=diamond*, mark size=2.5, mark options={solid,fill=darkslategray61}, only marks, forget plot]
table {%
	2.26666666666667 7
	2.26666666666667 24
	2.26666666666667 27
};
\addplot [semithick, darkslategray61, forget plot]
table {%
	-0.397333333333333 2605
	-0.136 2605
};
\addplot [semithick, darkslategray61, forget plot]
table {%
	-0.130666666666667 2200
	0.130666666666667 2200
};
\addplot [semithick, darkslategray61, forget plot]
table {%
	0.136 2164.5
	0.397333333333333 2164.5
};
\addplot [semithick, darkslategray61, forget plot]
table {%
	0.602666666666667 389
	0.864 389
};
\addplot [semithick, darkslategray61, forget plot]
table {%
	0.869333333333333 135
	1.13066666666667 135
};
\addplot [semithick, darkslategray61, forget plot]
table {%
	1.136 100
	1.39733333333333 100
};
\addplot [semithick, darkslategray61, forget plot]
table {%
	1.60266666666667 0
	1.864 0
};
\addplot [semithick, darkslategray61, forget plot]
table {%
	1.86933333333333 0
	2.13066666666667 0
};
\addplot [semithick, darkslategray61, forget plot]
table {%
	2.1360 0
	2.39733333333333 0
};
\end{axis} 

\end{tikzpicture}
        \caption{Requests not served by the fleet.}
        \label{fig:results:study:computational-study-unassigned-per-fleet}
     \end{subfigure}
     \hfill
     \begin{subfigure}[b]{0.48\textwidth}
        \centering
\begin{tikzpicture}

\definecolor{darkslategray38}{RGB}{38,38,38}
\definecolor{darkslategray61}{RGB}{61,61,61}
\definecolor{lightgray204}{RGB}{204,204,204}
\definecolor{peru22412844}{RGB}{224,128,44}
\definecolor{seagreen5814558}{RGB}{58,145,58}
\definecolor{steelblue49115161}{RGB}{49,115,161}

\begin{axis}[
width=0.8\textwidth,
height=8cm, 
axis line style={lightgray204},
legend cell align={left},
legend style={fill opacity=0.8, draw opacity=1, text opacity=1, at={(0.98,0.225)}, draw=lightgray204, nodes={scale=0.7, transform shape}},
tick align=outside,
x grid style={lightgray204},
xlabel=\textcolor{darkslategray38}{Fleet},
xmajorticks=true,
xmin=-0.5, xmax=2.5,
xtick style={color=darkslategray38},
xtick={0,1,2},
xticklabels={200,500,800},
y grid style={lightgray204},
ylabel=\textcolor{darkslategray38}{Total distance [km]},
ymajorgrids,
ymajorticks=true,
scaled ticks=false,
ymin=2665.68595, ymax=11786.93105,
ytick style={color=darkslategray38}
]
\path [draw=darkslategray61, fill=steelblue49115161, semithick]
(axis cs:-0.397333333333333,3230.9345)
--(axis cs:-0.136,3230.9345)
--(axis cs:-0.136,3311.788)
--(axis cs:-0.397333333333333,3311.788)
--(axis cs:-0.397333333333333,3230.9345)
--cycle;
\path [draw=darkslategray61, fill=peru22412844, semithick]
(axis cs:-0.130666666666667,3417.6735)
--(axis cs:0.130666666666667,3417.6735)
--(axis cs:0.130666666666667,3519.5405)
--(axis cs:-0.130666666666667,3519.5405)
--(axis cs:-0.130666666666667,3417.6735)
--cycle;
\path [draw=darkslategray61, fill=seagreen5814558, semithick]
(axis cs:0.136,3405.15125)
--(axis cs:0.397333333333333,3405.15125)
--(axis cs:0.397333333333333,3481.94625)
--(axis cs:0.136,3481.94625)
--(axis cs:0.136,3405.15125)
--cycle;
\path [draw=darkslategray61, fill=steelblue49115161, semithick]
(axis cs:0.602666666666667,7235.6475)
--(axis cs:0.864,7235.6475)
--(axis cs:0.864,7621.39375)
--(axis cs:0.602666666666667,7621.39375)
--(axis cs:0.602666666666667,7235.6475)
--cycle;
\path [draw=darkslategray61, fill=peru22412844, semithick]
(axis cs:0.869333333333333,7691.004)
--(axis cs:1.13066666666667,7691.004)
--(axis cs:1.13066666666667,7916.28175)
--(axis cs:0.869333333333333,7916.28175)
--(axis cs:0.869333333333333,7691.004)
--cycle;
\path [draw=darkslategray61, fill=seagreen5814558, semithick]
(axis cs:1.136,7582.47775)
--(axis cs:1.39733333333333,7582.47775)
--(axis cs:1.39733333333333,7878.8815)
--(axis cs:1.136,7878.8815)
--(axis cs:1.136,7582.47775)
--cycle;
\path [draw=darkslategray61, fill=steelblue49115161, semithick]
(axis cs:1.60266666666667,7498.8845)
--(axis cs:1.864,7498.8845)
--(axis cs:1.864,8791.38075)
--(axis cs:1.60266666666667,8791.38075)
--(axis cs:1.60266666666667,7498.8845)
--cycle;
\path [draw=darkslategray61, fill=peru22412844, semithick]
(axis cs:1.86933333333333,7656.65)
--(axis cs:2.13066666666667,7656.65)
--(axis cs:2.13066666666667,9080.59175)
--(axis cs:1.86933333333333,9080.59175)
--(axis cs:1.86933333333333,7656.65)
--cycle;
\path [draw=darkslategray61, fill=seagreen5814558, semithick]
(axis cs:2.136,7498.8845)
--(axis cs:2.39733333333333,7498.8845)
--(axis cs:2.39733333333333,8791.38075)
--(axis cs:2.136,8791.38075)
--(axis cs:2.136,7498.8845)
--cycle;
\draw[draw=darkslategray61,fill=steelblue49115161,line width=0.3pt] (axis cs:0,0) rectangle (axis cs:0,0);
\addlegendimage{ybar,ybar legend,draw=darkslategray61,fill=steelblue49115161,line width=0.3pt}
\addlegendentry{\mSolverSequentialMathHeur}

\draw[draw=darkslategray61,fill=peru22412844,line width=0.3pt] (axis cs:0,0) rectangle (axis cs:0,0);
\addlegendimage{ybar,ybar legend,draw=darkslategray61,fill=peru22412844,line width=0.3pt}
\addlegendentry{\mSolverDBILS}

\draw[draw=darkslategray61,fill=seagreen5814558,line width=0.3pt] (axis cs:0,0) rectangle (axis cs:0,0);
\addlegendimage{ybar,ybar legend,draw=darkslategray61,fill=seagreen5814558,line width=0.3pt}
\addlegendentry{\mSolverCombined}

\addplot [semithick, darkslategray61, forget plot]
table {%
	-0.266666666666667 3230.9345
	-0.266666666666667 3188.006
};
\addplot [semithick, darkslategray61, forget plot]
table {%
	-0.266666666666667 3311.788
	-0.266666666666667 3340.906
};
\addplot [semithick, darkslategray61, forget plot]
table {%
	-0.332 3188.006
	-0.201333333333333 3188.006
};
\addplot [semithick, darkslategray61, forget plot]
table {%
	-0.332 3340.906
	-0.201333333333333 3340.906
};
\addplot [semithick, darkslategray61, forget plot]
table {%
	0 3417.6735
	0 3328.976
};
\addplot [semithick, darkslategray61, forget plot]
table {%
	0 3519.5405
	0 3564.968
};
\addplot [semithick, darkslategray61, forget plot]
table {%
	-0.0653333333333333 3328.976
	0.0653333333333333 3328.976
};
\addplot [semithick, darkslategray61, forget plot]
table {%
	-0.0653333333333333 3564.968
	0.0653333333333333 3564.968
};
\addplot [semithick, darkslategray61, forget plot]
table {%
	0.266666666666667 3405.15125
	0.266666666666667 3325.22
};
\addplot [semithick, darkslategray61, forget plot]
table {%
	0.266666666666667 3481.94625
	0.266666666666667 3505.565
};
\addplot [semithick, darkslategray61, forget plot]
table {%
	0.201333333333333 3325.22
	0.332 3325.22
};
\addplot [semithick, darkslategray61, forget plot]
table {%
	0.201333333333333 3505.565
	0.332 3505.565
};
\addplot [semithick, darkslategray61, forget plot]
table {%
	0.733333333333333 7235.6475
	0.733333333333333 6942.544
};
\addplot [semithick, darkslategray61, forget plot]
table {%
	0.733333333333333 7621.39375
	0.733333333333333 7863.583
};
\addplot [semithick, darkslategray61, forget plot]
table {%
	0.668 6942.544
	0.798666666666667 6942.544
};
\addplot [semithick, darkslategray61, forget plot]
table {%
	0.668 7863.583
	0.798666666666667 7863.583
};
\addplot [semithick, darkslategray61, forget plot]
table {%
	1 7691.004
	1 7516.292
};
\addplot [semithick, darkslategray61, forget plot]
table {%
	1 7916.28175
	1 8039.102
};
\addplot [semithick, darkslategray61, forget plot]
table {%
	0.934666666666667 7516.292
	1.06533333333333 7516.292
};
\addplot [semithick, darkslategray61, forget plot]
table {%
	0.934666666666667 8039.102
	1.06533333333333 8039.102
};
\addplot [black, mark=diamond*, mark size=2.5, mark options={solid,fill=darkslategray61}, only marks, forget plot]
table {%
	1 7274.998
	1 7132.082
	1 7317.461
	1 7082.943
};
\addplot [semithick, darkslategray61, forget plot]
table {%
	1.26666666666667 7582.47775
	1.26666666666667 7256.416
};
\addplot [semithick, darkslategray61, forget plot]
table {%
	1.26666666666667 7878.8815
	1.26666666666667 7977.454
};
\addplot [semithick, darkslategray61, forget plot]
table {%
	1.20133333333333 7256.416
	1.332 7256.416
};
\addplot [semithick, darkslategray61, forget plot]
table {%
	1.20133333333333 7977.454
	1.332 7977.454
};
\addplot [black, mark=diamond*, mark size=2.5, mark options={solid,fill=darkslategray61}, only marks, forget plot]
table {%
	1.26666666666667 7102.662
	1.26666666666667 7136.564
	1.26666666666667 7056.376
};
\addplot [semithick, darkslategray61, forget plot]
table {%
	1.73333333333333 7498.8845
	1.73333333333333 6990.145
};
\addplot [semithick, darkslategray61, forget plot]
table {%
	1.73333333333333 8791.38075
	1.73333333333333 10677.168
};
\addplot [semithick, darkslategray61, forget plot]
table {%
	1.668 6990.145
	1.79866666666667 6990.145
};
\addplot [semithick, darkslategray61, forget plot]
table {%
	1.668 10677.168
	1.79866666666667 10677.168
};
\addplot [black, mark=diamond*, mark size=2.5, mark options={solid,fill=darkslategray61}, only marks, forget plot]
table {%
	1.73333333333333 10854.889
	1.73333333333333 10988.094
};
\addplot [semithick, darkslategray61, forget plot]
table {%
	2 7656.65
	2 7066.921
};
\addplot [semithick, darkslategray61, forget plot]
table {%
	2 9080.59175
	2 11114.755
};
\addplot [semithick, darkslategray61, forget plot]
table {%
	1.93466666666667 7066.921
	2.06533333333333 7066.921
};
\addplot [semithick, darkslategray61, forget plot]
table {%
	1.93466666666667 11114.755
	2.06533333333333 11114.755
};
\addplot [black, mark=diamond*, mark size=2.5, mark options={solid,fill=darkslategray61}, only marks, forget plot]
table {%
	2 11372.329
};
\addplot [semithick, darkslategray61, forget plot]
table {%
	2.26666666666667 7498.8845
	2.26666666666667 6990.145
};
\addplot [semithick, darkslategray61, forget plot]
table {%
	2.26666666666667 8791.38075
	2.26666666666667 10717.379
};
\addplot [semithick, darkslategray61, forget plot]
table {%
	2.20133333333333 6990.145
	2.332 6990.145
};
\addplot [semithick, darkslategray61, forget plot]
table {%
	2.20133333333333 10717.379
	2.332 10717.379
};
\addplot [black, mark=diamond*, mark size=2.5, mark options={solid,fill=darkslategray61}, only marks, forget plot]
table {%
	2.26666666666667 11011.782
	2.26666666666667 11255.444
};
\addplot [semithick, darkslategray61, forget plot]
table {%
	-0.397333333333333 3250.2875
	-0.136 3250.2875
};
\addplot [semithick, darkslategray61, forget plot]
table {%
	-0.130666666666667 3480.587
	0.130666666666667 3480.587
};
\addplot [semithick, darkslategray61, forget plot]
table {%
	0.136 3454.3225
	0.397333333333333 3454.3225
};
\addplot [semithick, darkslategray61, forget plot]
table {%
	0.602666666666667 7386.651
	0.864 7386.651
};
\addplot [semithick, darkslategray61, forget plot]
table {%
	0.869333333333333 7844.929
	1.13066666666667 7844.929
};
\addplot [semithick, darkslategray61, forget plot]
table {%
	1.136 7755.6905
	1.39733333333333 7755.6905
};
\addplot [semithick, darkslategray61, forget plot]
table {%
	1.60266666666667 7841.013
	1.864 7841.013
};
\addplot [semithick, darkslategray61, forget plot]
table {%
	1.86933333333333 8028.259
	2.13066666666667 8028.259
};
\addplot [semithick, darkslategray61, forget plot]
table {%
	2.136 7835.80150
	2.39733333333333 7835.8015
};
\end{axis}
 
\end{tikzpicture}
        \caption{Total distance traveled by the fleet.}
        \label{fig:results:study:computational-study-distance-per-fleet}
     \end{subfigure}
     \caption{Computational study final results per fleet size.}
     \label{fig:results:study:computational-study-per-fleet}
\end{figure}
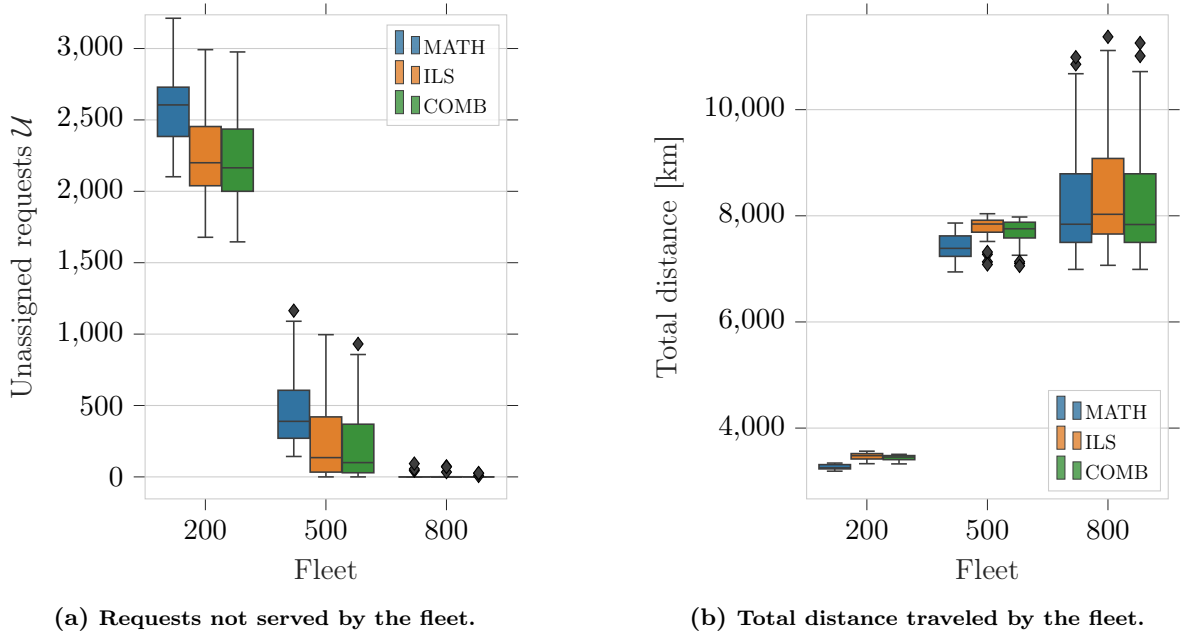

Figure~\ref{fig:results:study:computational-study-per-fleet} shows the performance of $\mSolverSequentialMathHeur$, $\mSolverDBILS$, and $\mSolverCombined$ over the set of mid-sized instances $\textsc{Mid}^{25}$ with 25\% of the total requests for different fleet sizes. As can be seen, $\mSolverDBILS$ and $\mSolverCombined$ outperform $\mSolverSequentialMathHeur$ on our primary objective, which is to minimize the number of unassigned requests, see Figure~8a. Figure~8b shows our secondary objective which is to minimize the total traveled distance. Here, it is not possible to provide a direct comparison as the total traveled distance could only be compared for the same primary objective values. However, it is worth notifying that $\mSolverCombined$, which performs equal or better than $\mSolverDBILS$ on our primary objective, also outperforms $\mSolverDBILS$ on the secondary objective. 

Beyond these objective-focused observations, we identify two additional effects: in settings with an abundant supply of vehicles, all algorithms yield a similar performance on the primary objective. As soon as vehicle availability becomes a bottleneck, we observe a significant performance difference between the studied algorithms. Analyzing Figure~8a, we further observe that the share of unassigned requests remains too high from a practical perspective as an operator would aim at expanding its fleet size to provide a better service level. With these observations in mind, we focus the remaining discussions on the setting with a fleet size of 500 vehicles as i) the setting with 800 vehicles does not allow for meaningful algorithmic studies as it compensates bad algorithmic decisions by abundant vehicle supply, and ii) the setting with 200 vehicles leads to unrealistically high rates of unserved requests from a practitioners perspective.

Figure~\ref{fig:results:study:computational-study-progress-unassigned} and \ref{fig:results:study:computational-study-progress-distance} show the search progress over the runtime for each buffer setting $\mBuffer=1,\ldots,6$. The result of $\mSolverSequentialMathHeur$, which finishes in a few minutes, is shown as a horizontal line for comparison. The combined approach, $\mSolverCombined$, is warm-started with the solution from $\mSolverSequentialMathHeur$. In terms of unassigned requests, we observe that the integrated metaheuristic $\mSolverDBILS$ reaches and improves the solution of $\mSolverSequentialMathHeur$ in at most 10 minutes on average. $\mSolverCombined$ is able to improve on the sequential solution and produces better results than $\mSolverDBILS$ approach in the runtime of one hour. However, the difference is getting small with higher buffer settings $\mBuffer \geq 3$. In any case, $\mSolverCombined$ converges faster to its final value then $\mSolverDBILS$.

Focusing on the total distance traveled, we observe an increase during the search in $\mSolverCombined$ when started with the solution from $\mSolverSequentialMathHeur$. This effect correlates with a decrease in unassigned requests; however, we see stagnation after a few minutes, i.e., even though improvements in terms of unassigned requests are found, the total distance traveled by the fleet does not increase significantly. This implies that the metaheuristic is able to find more efficient routing decisions, which is especially visible for higher buffer settings, e.g., $\mBuffer\geq 5$, where the total distance almost reaches initial solution found by $\mSolverSequentialMathHeur$. Remarkably, the $\mSolverCombined$ approach succeeds in serving more customers at a lower totaled traveled distance compared to the $\mSolverDBILS$ approach.

\begin{figure}[!ht]
    \centering
    \input{figures/computational_study_progress_unassigned_2x3}
    \caption{Search progress (best solution) over time regarding unassigned requests per buffer setting for 500 vehicles.}
    \label{fig:results:study:computational-study-progress-unassigned}
\end{figure}

\begin{figure}[!ht]
    \centering
    \input{figures/computational_study_progress_distance_2x3}
    \caption{Search progress (best solution) over time regarding total distance traveled per buffer setting for 500 vehicles.}
    \label{fig:results:study:computational-study-progress-distance}
\end{figure}

\begin{figure}
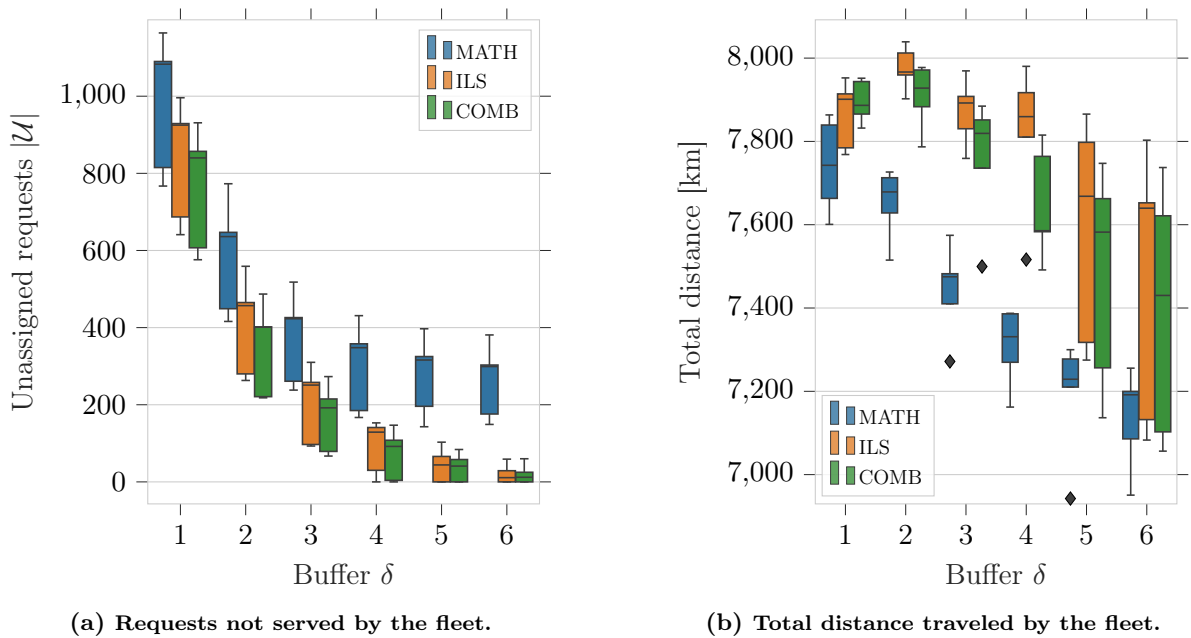

    \centering
     \begin{subfigure}[b]{0.48\textwidth}
        \centering
        \input{figures/computational_study_unassigned_per_buffer}
        \caption{Requests not served by the fleet.}
        \label{fig:results:study:computational-study-unassigned-per-buffer}
     \end{subfigure}
     \hfill
     \begin{subfigure}[b]{0.48\textwidth}
        \centering
        \input{figures/computational_study_distance_per_buffer}
        \caption{Total distance traveled by the fleet.}
        \label{fig:results:study:computational-study-distance-per-buffer}
     \end{subfigure}
     \caption{Computational study final results per buffer length.}
     \label{fig:results:study:computational-study-per-buffer}
\end{figure}

Figure~\ref{fig:results:study:computational-study-per-buffer} shows the final solutions for the instances with 500 vehicles, separated by the buffer length $\mBuffer$. We see a similar performance of $\mSolverDBILS$ and $\mSolverCombined$ for unassigned requests. $\mSolverSequentialMathHeur$, however, seems to stagnate with higher buffer settings. In the case of total distance traveled by the fleet, we can see that $\mSolverSequentialMathHeur$ produces solutions with less distance traveled. This is not surprising, as these solutions serve less requests than the solutions found by the metaheuristic. Furthermore, we can observe that the best performing $\mSolverCombined$ approach does tend to produce solutions with not only more requests served but less total distance traveled by the fleet.

In summary, we see $\mSolverSequentialMathHeur$ produces good solutions for scenarios with large fleets, where most of the requests can be served. In these cases, improving these solutions using our metaheuristic ($\mSolverCombined$) may not result in significant improvements, especially with higher buffer settings (for more details, see Appendix~\ref{appendix:pooling}). For the under-saturated case, where the fleet is sparse, we can see that solutions from $\mSolverSequentialMathHeur$ are outperformed by $\mSolverDBILS$. 
Comparing the results of the metaheuristic $\mSolverDBILS$ with the default initial solution generation and warm-started with the solution by $\mSolverSequentialMathHeur$, i.e., $\mSolverCombined$, we can see that the approach benefit from a good initial solution, resulting in better final solutions -- 6.8\% fewer unassigned requests on average -- after a time limit of just one hour.

\section{PDPTW benchmark}\label{appendix:pdptw-benchmark}

We evaluate the efficacy of our approach by comparing its performance solving \gls{pdptw} instances. \cite{SartoriBuriol2020} recently created a new benchmark based on real-world street network data of Barcelona, Berlin, New York City, and Porto Alegre, ranging from 100 to 5000 nodes (50 to 2500 requests). Each consists of 25 instances, 5-7 per city. Current best-known solutions are curated by the benchmark creator and available online at \url{https://github.com/cssartori/pdptw-instances}. 

The classic \gls{pdptw} uses a hierarchical objective, where we prioritize minimizing the number of vehicles used before the total distance traveled. Furthermore, no unassigned requests are allowed in a feasible solution. Herein, we include a fleet minimization component to our \gls{ils}-based metaheuristic approach. We refer to Appendix~\ref{appendix:fleet-min-algorithm} for details.

Note that we limit our discussions to the benchmark of \cite{SartoriBuriol2020}. In a preliminary analysis, we looked at the structure of the best solutions produced in the literature for the \cite{LiLim2003} and \cite{SartoriBuriol2020} benchmarks. We assessed their similarity to our large-scale ride-sharing case (see Appendix~\ref{appendix:solution-structure}). Our findings show that the classical \cite{LiLim2003} instances tend to have a higher tendency to result in routing decisions with long and sometimes route-spanning blocks. In contrast, solutions for the \cite{SartoriBuriol2020} benchmark exhibit a higher rate of blocks assigned to routes. This property overlaps with the solutions resulting from our large-scale ride-sharing instances. To this end, we focused on the \cite{SartoriBuriol2020} benchmark instances for our comparative analysis.

Table \ref{tab:sartoriburiol} compares our results with the literature for the \cite{SartoriBuriol2020} benchmark. As is usually done in the \gls{pdptw} literature, we report our results as accumulated values (the number of vehicles and total cost, i.e., distance traveled) per instance size. Note that the best-known solutions reported on \url{https://github.com/cssartori/pdptw-instances} include unpublished work without any additional information regarding runtime and hardware, and a recent contribution by \cite{VadsethAnderssonEtAl2023}, which improved 203 instances by warm-starting their approach using the previously best-known solution. To our knowledge, the only published results produced without using best-known solutions are from \cite{SartoriBuriol2020}.

Our results show the competitiveness of our approach, especially for larger instance sizes, where we improved on the previous results by 3.1\% in terms of the cumulative number of vehicles. Furthermore, we contribute to the research by providing 107 new best solutions, with 90 reducing the number of vehicles.

\begin{table}[tbh!]

\newcommand{\emphBetterThanSB}[1]{\emphBold{#1}}
\newcommand{\emphBetterThanBKS}[1]{\emphBoldUnderlined{#1}}
\begin{tabular}{lrr|rrrrrrrr}
\toprule
& \multicolumn{2}{c}{BKS} & \multicolumn{2}{c}{SB*} & \multicolumn{2}{c}{HS*} & \multicolumn{2}{c}{HS**} \\ \cmidrule(l){2-3} \cmidrule(l){4-5} \cmidrule(l){6-7} \cmidrule(l){8-9} 
\multicolumn{1}{c}{Instance} & \multicolumn{1}{c}{Veh}     & \multicolumn{1}{c}{Cost} & \multicolumn{1}{c}{Veh}     & \multicolumn{1}{c}{Cost}         & \multicolumn{1}{c}{Veh} & \multicolumn{1}{c}{Cost} & \multicolumn{1}{c}{Veh} & \multicolumn{1}{c}{Cost} \\ \midrule
n100  & 164   & 25.264 & 164   & 25.388 & 164   & \emphBetterThanSB{25.349} & 164 & \emphBetterThanSB{25.326} \\
n200  & 332   & 46.617 & 337   & 46.587 & \emphBetterThanSB{332}   & 46.846 & \emphBetterThanSB{332}   & 46.686 \\
n400  & 585   & 85.277 & 589   & 85.887 & 589   & 87.604 & \emphBetterThanSB{587}   & 86.834 \\
n600  & 828   & 122.156 & 840   & 122.130 & \emphBetterThanSB{838}   & 126.081 & \emphBetterThanSB{834}   & 124.996 \\
n800  & 1134  & 163.806 & 1150  & 163.341 & \emphBetterThanSB{1147}  & 170.952 & \emphBetterThanSB{1147}  & 169.756 \\
n1000 & 1383  & 222.846 & 1401  & 226.228 & \emphBetterThanSB{1389}  & 232.096 & \emphBetterThanSB{1385}  & 230.305 \\
n1500 & 2078  & 297.367 & 2115  & 303.478 & \emphBetterThanSB{2082}  & 313.057 & \emphBetterThanBKS{2076}  & 311.215 \\
n2000 & 2861  & 412.036 & 2924  & 425.343 & \emphBetterThanSB{2902}  & 438.839 & \emphBetterThanSB{2886}  & 433.970 \\
n2500 & 3135  & 493.731 & 3201  & 506.436 & \emphBetterThanSB{3149}  & 519.504 & \emphBetterThanBKS{3128}  & 513.210 \\
n3000 & 4137  & 583.567 & 4276  & 625.072 & \emphBetterThanBKS{4130}  & \emphBetterThanSB{614.621} & \emphBetterThanBKS{4105}  & 605.405 \\
n4000 & 5534  & 761.630 & 5944  & 866.430 & \emphBetterThanSB{5605}  & \emphBetterThanSB{820.816} & \emphBetterThanSB{5541}  & \emphBetterThanSB{796.389} \\
n5000 & 6262  & 918.989 & 6802  & 1.089.677 & \emphBetterThanSB{6443}  & \emphBetterThanSB{1.029.592} & \emphBetterThanSB{6335}  & \emphBetterThanSB{982.320} \\ \bottomrule
\multicolumn{9}{l}{\scriptsize BKS: best-known solutions from \url{https://github.com/cssartori/pdptw-instances}} \\
\multicolumn{9}{l}{ \scriptsize(last-access: 2024-04-26); SB: \cite{SartoriBuriol2020}; HS: this work} \\
\multicolumn{9}{l}{ \scriptsize * time limit of 5, 15, 15, 30, and 60 minutes for n100, n200, n400, n600, and n800--n5000} \\
\multicolumn{9}{l}{ \scriptsize ** extended time limit of 15, 45, 45, 90, and 180 minutes for n100, n200, n400, n600, and n800--n5000} \\
\end{tabular}

\caption{\cite{SartoriBuriol2020} benchmark results (cumulative values for the number of vehicles and cost). We mark results in bold if they improve the results reported by \cite{SartoriBuriol2020}, and underline results improving on the previously best-known solutions.}
\label{tab:sartoriburiol}
\end{table}

\subsection{Managerial Analyses}
In the following, we use our $\mSolverDBILS$ to study the large-scale instances introduced in Section~\ref{subsection:case-study} from a managerial perspective. In this context, we discuss the impact of fleet-sizing and delays incurred on the customer side.

\begin{figure}[tbh!]
    \centering
    \input{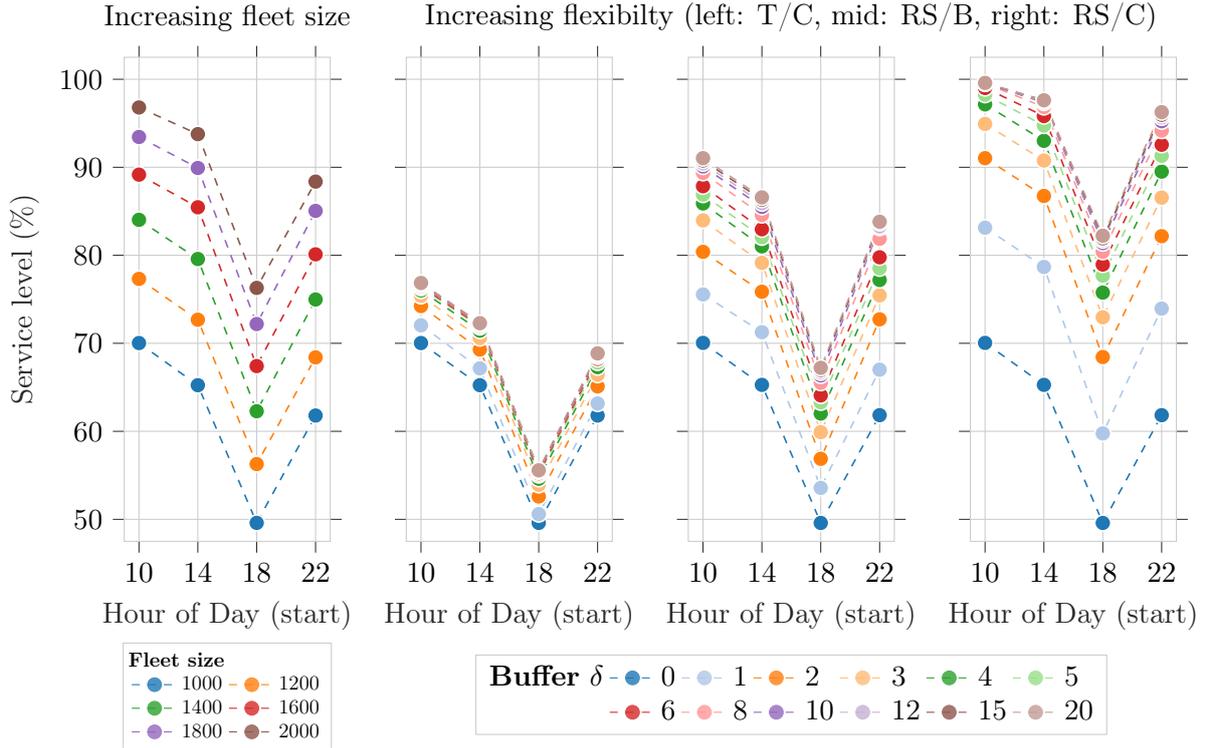}
    \caption{Change in service level comparison between fleet size and flexibility (T (taxi): single requests; RS (ride-sharing): up to three requests / B: fixed pickup, variable dropoff; C~ variable pickup and dropoff).}
    \label{fig:results:nyc:comparison_fleet_buffer_weekdays}
\end{figure}

\paragraph{Comparison of fleet variations to ride-sharing flexibility.}
We evaluate the benefits of introducing flexibility for the service level by comparing our results to a commonly used -- although expensive -- alternative: increasing the fleet size. We focus on a single day and consider fleets of size 1000 to 2000 vehicles, with 200-vehicle increments. Figure~\ref{fig:results:nyc:comparison_fleet_buffer_weekdays} shows the corresponding results for four scenarios: a) the aforementioned fleet increase of the classic taxi service, b) taxi operations with possible delayed pickup (and delivery), c) ride-sharing with a fixed pickup time, but variable delivery, and d) ride-sharing with flexible pickup and delivery time. We observe a slower growth in service level when the fleet size increases, whereas increased allowed delay introduces more considerable gains with lower values and diminishing increases with larger values. 
For the most flexible option (ride-sharing with variable pickup and delivery time up to $\delta$ minutes), we see that even a small buffer of two minutes allows a similar service level than an increase of the fleet by 600 additional vehicles for the time frame starting at 6 p.m. (18h). 

\paragraph{Buffer effects on ride-sharing operations.}

We now continue with the most flexible case, i.e., ride-sharing with increasingly flexible pickup and delivery times. Figure~\ref{fig:results:nyc:service_level_per_buffer_weekdays} shows the service level -- the percentage of requests served by the fleet -- for increasing buffer levels $\mBuffer$, i.e., increasing flexibility of arrival delays. Given the fleet limitation of 1000 vehicles, not all requests can be served for dense time-frames, e.g., 6 p.m. (18h). We can see a larger relative gain in service level with smaller allowed delays (buffer). Even a two-minute delay allows for an increase from approximately 60 \% to 80\%, as shown for 10 a.m. (10h) and 2 p.m. (14h). 
We further observe a larger spread in the service level for the latest time frame at 10 p.m. (22h). This effect originates primarily from the results of the Thursday and Friday instances, which tend to compose of more a higher frequency of requests compared to Monday to Wednesday. 

\begin{figure}[tbh!]
    \centering
    \input{figures/service_level_per_buffer_weekdays}
    \caption{Service level changes with increasing allowed delay (buffer) 
for weekdays per hour of the day with a fixed fleet of 1000 vehicles.}
    \label{fig:results:nyc:service_level_per_buffer_weekdays}
\end{figure}

Figure~\ref{fig:results:nyc:distance_traveled_per_request_served_per_buffer_weekdays} shows the distance traveled per request served to increase the flexibility of arrival delays. We can observe, as expected, that in the lower density time-frames at 10 a.m. (10h) and 2 p.m. (14h), the distance traveled per request served is larger than for the dense evening time frame at 6 p.m. Furthermore, we see an increase in average distance traveled when introducing flexibility; however, this effect does not necessarily continue but stagnates after the 2-minute allowed delay (see 18h) or even decreases for the early time frame (10h). Comparing the distance progression with the service level changes in Figure~\ref{fig:results:nyc:service_level_per_buffer_weekdays}, we can see no clear correlation between those indicators, especially for the dense time frame, with roughly similar average distance values of 800 meters for 2 to 20-minute flexibility, with an increase of 10 \% of service level for the same range. 

\begin{figure}[tbh!]
    \centering
    \input{figures/distance_traveled_per_fulfilled_requests_weekdays}
    \caption{Distance traveled per request served with increasing allowed delay (buffer) 
for weekdays per hour of the day with a fixed fleet of 1000 vehicles.}
    \label{fig:results:nyc:distance_traveled_per_request_served_per_buffer_weekdays}
\end{figure}

Figure~\ref{fig:results:nyc:average_delay_per_buffer_weekdays} shows the effects of introducing ride-sharing and the induced allowed delay on the experienced delay of customers. The average experienced delay is always below half the permitted delay and settles at around 6 minutes for the highest delay setting of 20 minutes. For the high-density time frames at 6 p.m., we see a lower increase in experience delay with higher permitted delay settings.

\begin{figure}[tbh!]
    \centering
    \input{figures/average_delay_per_buffer_weekdays}
    \caption{Average delay per customer due to introduced flexibility (buffer) 
 for weekdays and grouped by hour of day with a fixed fleet of 1000 vehicles.}
    \label{fig:results:nyc:average_delay_per_buffer_weekdays}
\end{figure}

\FloatBarrier

\section{Conclusion}\label{section:conclusion}
In this paper, we presented an algorithmic framework to solve very large-scale instances of an \gls{pdptw} variant arising in the context of urban ride-sharing systems. Herein, we developed three approaches: i) a decomposition-based matheuristic, which can solve instances with up to 5000 requests in a few minutes; ii) a \gls{ils}-based metaheuristic to handle very large-scale instances with more than 20 thousand requests; and iii) a hybrid approach, where we warm-start the metaheuristic with the solution from the matheuristic.

We conducted a thorough computational study using this algorithmic framework to derive insights into the characteristics of the proposed algorithms. We saw that the metaheuristic improved on the matheuristic within a time limit of 15 minutes. The warm-started hybrid approach produced better results, with up to 6.8\% fewer unassigned requests on average. 
Next, we compared our metaheuristic approach to results for the \gls{pdptw} benchmark data set comprising instances with up to 2500 requests. The results showed the algorithm's competitiveness and found 107 new best solutions, with 90 solutions reducing the number of vehicles. 
Finally, we applied our algorithm to very large-scale instances with up to 21375 requests. Our results show that by requesting passengers to allow flexibility in their delivery by just 2 minutes, the mobility provider can operate on a comparative service level as with a 50\% larger fleet. Additionally, on average, the passengers would only delay their arrival by less than a minute. 

Future research may expand on the instances' scale by using only a subset of trips but expanding the time frame to multiple hours or a full day. Herein, new challenges for the decomposition approach may arise regarding split and recombination policy.

\section*{Acknowledgements}{
	This research has been funded by the Deutsche Forschungsgemeinschaft (DFG, German Research Foundation) - Project number 449261765 (BalSAM). This support is gratefully acknowledged.
}


%
\singlespacing{
\bibliographystyle{model5-names}
\bibliography{./bibliography/main,./bibliography/oscm}} 
\newpage
\onehalfspacing
\begin{appendices}
	\normalsize
	\section{Computational study for the hypergraph matching}
\label{appendix:pooling}

In Section~\ref{section:methodology:matheuristic}, we described the best-performing hypergraph matching approach we found in our research. In the following, we provide more details on the variations we tested, the experiments conducted, and the results obtained. 
We start by introducing a simple greedy procedure for selecting a matching.
Next, we provide a \gls{wsp} formulation based on the \gls{wsp} formulation from Section~\ref{fig:sequential-pooling-and-dispatching:pooling}. 
Next, we show the alternative fractional selection policy, based on randomized rounding. Finally, we conclude this section by providing detailed results on the performance of these approaches, which justifies the final selection for our matheuristic approach.

\subsection{Greedy matching}
We initially considered a naive approach to handle the matching process. Herein, we first simply sort all hyperedges by one of the weight functions. Then, we traverse this list and pick the next best item if none of the associated requests have been part of a previously selected hyperedge. Finally, we stop if all requests are part of selected hyperedges, or we reach the end of the list.

\subsection{Weighted set partitioning problem}

A matching $\mMatching$ for hypergraph $\mHyperGraph$ is defined as a subset of hyperedges $\mMatching \subset \mSetHyperEdges$ where requests are only part of at most one hyperedge. To find $\mMatching$, we tested different approaches. Apart from the greedy selection approach, we considered two two-step approaches. One is the formulation as a \gls{wsc} problem, as described in Section~\ref{section:methodology:matheuristic:matching}. In the other approach, we solve a continuous \gls{wsp} problem to select a subset of hyperedges that may be connected to the same request nodes, i.e., a fractional matching $\overline{\mMatching}$. In a second step, we derive matching $\mMatching$ based on the solution of the \gls{wsp} using either a greedy selection procedure, or a randomized rounding approach. 

We can formulate the \gls{wsp} by replacing \eqref{eq:cover} from Section~\ref{section:methodology:matheuristic:matching} with
\begin{align}
&\sum_{\mHyperEdge{} \in \mSetHyperEdges'} a_{\mRequest{}\mHyperEdge{}} x_{\mHyperEdge{}} = 1 &&\forall \mRequest{} \in \mSetRequests \label{eq:wsp:cover}
\end{align}

Note that in preliminary experiments we tested an ILP formulation, to avoid fractional solutions. However, run times of over one hour turned out to be impractical, and the solutions obtained did not considerable improve on the two-step approach.

\paragraph{Selection from fractional solution.}

When solving the relaxed \gls{wsc} (\gls{wsp}) problem, we need to deal with fractional solution values. We test two strategies: a) a simple greedy selection procedure (see Section~\ref{section:methodology:matheuristic} and b) a standard randomized rounding scheme.

The randomized rounding scheme follows the standard technique introduced by \cite{RaghavanThomson1987}. Herein, a cover (i.e., a pooling option) is selected at random with a probability defined by its fractional solution value. The resulting selection may cover requests multiple times. We enforce non-overlapping pooling options by using a simple repair procedure, where the request will be removed from all but one of the selected options.

\begin{table}[tbh!]
\small
\begin{tabular}{cccccc}
\toprule
Components & $\mMatchingGreedy$ & $\mMatchingWSCGreedy$ & $\mMatchingWSCSRR$ & $\mMatchingWSPGreedy$ & $\mMatchingWSPSRR$ \\
\midrule
Greedy matching & x & & & & \\
Set Covering & & x & x & & \\
Set Partitioning & & & & x & x \\
Greedy selection & & x & & x & \\
Randomized rounding & & & x & & x \\
\bottomrule
\end{tabular}
\caption{Summary of the proposed pooling selection procedures.}
\label{tab:matching-procedures}
\end{table}

\subsection{Results}

We test our approach on mid-sized instances, considering 25\% of trips, for the high-volume time frame of 18-19h. Our instances consist of approximately 5000 requests. Note that we preprocess the request pooling options. This process requires several hours of computation, espacially for larger settings of $\mBuffer$. For this reason, we limited $\mBuffer \leq 6$. Furthermore, we do not report any results on a setting of $\mBuffer = 0$, as the prepocessing showed that no pooling option with two or more requests could be found.

In Table~\ref{tab:matching-cmp} we show the average results over all instances. These results indicated that the naive greedy selection procedure outperforms the \gls{wsc} and \gls{wsp} approaches.
Finally, Table~\ref{tab:matching-cmp-fleet-200}--Table~\ref{tab:matching-cmp-fleet-800} show the results for each fleet size setting separately. 

\begin{table}[tbh!]
\small
\begin{tabular}{llrrrrr}
\toprule
$\omega$ & $\rho$ & $\mMatchingGreedy$ & $\mMatchingWSCGreedy$ & $\mMatchingWSCSRR$ & $\mMatchingWSPGreedy$ & $\mMatchingWSPSRR$ \\
\midrule
\multicolumn{2}{l}{$\omega^1$} & 1180.0/7194.355 & 1334.6/7385.635 & 1438.9/7467.493 & 1311.7/7341.306 & 1373.3/7466.447 \\
\multirow{11}{*}{$\omega^4(\rho)$} & 0.0 & 870.7/6838.231 & 890.7/6703.879 & 1018.5/7017.773 & 895.4/6708.495 & 1021.9/7040.573 \\
& 0.1 & 871.4/6838.329 & 888.4/6701.300 & 1016.6/7012.572 & 894.0/6701.676 & 1014.9/7043.129 \\
& 0.2 & 871.6/6836.022 & 887.2/6699.581 & 1014.8/7012.429 & 891.9/6699.279 & 1013.0/7031.272 \\
& 0.3 & 873.1/6842.598 & 885.6/6701.048 & 1010.2/7009.539 & 890.6/6701.106 & 1008.7/7036.593 \\
& 0.4 & 875.7/6845.000 & 882.9/6698.400 & 1003.1/7004.935 & 888.0/6699.360 & 1005.3/7041.258 \\
& 0.5 & 883.1/6854.974 & 881.5/6701.853 & 1001.5/7009.856 & 886.8/6697.925 & 1001.9/7026.194 \\
& 0.6 & 891.4/6878.164 & 879.6/6710.580 & 999.1/6996.288 & 883.8/6709.140 & 994.2/7028.833 \\
& 0.7 & 908.0/6916.147 & 880.1/6713.553 & 993.5/6995.961 & 884.2/6715.881 & 995.3/7024.404 \\
& 0.8 & 942.4/6968.142 & 881.7/6737.104 & 988.1/6990.832 & 886.1/6733.913 & 986.0/7020.154 \\
& 0.9 & 1017.0/7097.965 & 891.5/6778.400 & 974.4/6984.691 & 890.2/6777.310 & 973.0/6992.534 \\
& 1.0 & 1153.5/7254.055 & 1180.3/7196.074 & 1675.3/7336.007 & 1180.0/7194.355 & 1675.5/7335.932 \\
\bottomrule
\end{tabular}
\caption{Comparison of the matching approaches averaged over all fleet sizes ($200,\ldots,800$) buffer settings $\mBuffer=1,\ldots,6$ ($\mUnassigned$/Cost).}
\label{tab:matching-cmp}
\end{table}

\begin{table}[tbh!]
\small
\begin{tabular}{llrrrrr}
\toprule
$\omega$ & $\rho$ & $\mMatchingGreedy$ & $\mMatchingWSCGreedy$ & $\mMatchingWSCSRR$ & $\mMatchingWSPGreedy$ & $\mMatchingWSPSRR$ \\
\midrule
\multicolumn{2}{l}{$\omega^1$} & 2985.5/3278.944 & 3069.4/3322.791 & 3090.8/3361.338 & 3039.9/3322.971 & 3041.1/3351.332 \\
\multirow{11}{*}{$\omega^4(\rho)$} & 0.0 & 2489.9/3195.219 & 2620.8/3254.686 & 2685.7/3300.255 & 2648.5/3260.037 & 2717.7/3299.972 \\
& 0.1 & 2490.0/3195.425 & 2618.0/3254.441 & 2684.3/3300.273 & 2647.2/3253.786 & 2714.2/3294.868 \\
& 0.2 & 2490.4/3196.507 & 2618.0/3250.695 & 2683.0/3295.981 & 2644.9/3254.973 & 2710.3/3297.761 \\
& 0.3 & 2492.7/3198.575 & 2615.7/3255.910 & 2682.9/3296.203 & 2642.3/3257.690 & 2708.6/3304.740 \\
& 0.4 & 2497.5/3198.114 & 2610.7/3262.320 & 2675.2/3309.270 & 2640.8/3261.455 & 2703.8/3311.012 \\
& 0.5 & 2505.7/3208.025 & 2608.4/3265.880 & 2676.0/3305.797 & 2638.6/3266.643 & 2702.4/3304.229 \\
& 0.6 & 2522.3/3209.043 & 2608.3/3273.450 & 2674.5/3310.582 & 2636.7/3271.027 & 2696.7/3310.412 \\
& 0.7 & 2555.1/3219.217 & 2611.4/3276.985 & 2676.4/3310.786 & 2639.7/3271.269 & 2703.4/3310.447 \\
& 0.8 & 2605.6/3225.229 & 2616.3/3289.606 & 2676.7/3314.049 & 2639.7/3290.767 & 2698.2/3318.870 \\
& 0.9 & 2714.2/3234.905 & 2630.3/3298.426 & 2675.4/3312.252 & 2645.4/3298.910 & 2691.4/3301.927 \\
& 1.0 & 2914.1/3240.492 & 2970.8/3281.692 & 3214.6/3111.353 & 2965.2/3284.368 & 3213.6/3111.950 \\
\bottomrule
\end{tabular}
\caption{Comparison of the matching approaches for the 200 vehicle instances, averaged over buffer settings $\mBuffer=1,..,6$.}
\label{tab:matching-cmp-fleet-200}
\end{table}

\begin{table}[tbh!]
\small
\begin{tabular}{llrrrrr}
\toprule
$\omega$ & $\rho$ & $\mMatchingGreedy$ & $\mMatchingWSCGreedy$ & $\mMatchingWSCSRR$ & $\mMatchingWSPGreedy$ & $\mMatchingWSPSRR$ \\
\midrule
\multicolumn{2}{l}{$\omega^1$} & 2186.1/4818.741 & 2323.8/4870.416 & 2368.3/4922.043 & 2293.9/4867.301 & 2309.1/4917.445 \\
\multirow{11}{*}{$\omega^4(\rho)$} & 0.0 & 1640.0/4665.494 & 1749.3/4759.670 & 1856.7/4833.735 & 1777.8/4761.643 & 1887.0/4831.570 \\
& 0.1 & 1640.7/4669.768 & 1745.1/4760.180 & 1855.1/4828.085 & 1775.4/4758.757 & 1879.9/4837.343 \\
& 0.2 & 1642.3/4662.946 & 1744.3/4759.608 & 1849.6/4827.459 & 1772.2/4762.537 & 1875.3/4830.555 \\
& 0.3 & 1643.1/4671.965 & 1744.5/4756.601 & 1847.8/4837.390 & 1771.1/4758.564 & 1875.2/4839.483 \\
& 0.4 & 1649.1/4672.658 & 1738.6/4763.288 & 1840.5/4838.725 & 1768.1/4764.814 & 1870.4/4840.993 \\
& 0.5 & 1660.9/4682.427 & 1736.5/4767.726 & 1839.8/4842.356 & 1767.2/4767.261 & 1863.7/4843.515 \\
& 0.6 & 1678.0/4690.973 & 1735.8/4778.772 & 1837.8/4848.617 & 1763.5/4776.552 & 1858.2/4854.202 \\
& 0.7 & 1708.5/4717.316 & 1737.7/4786.482 & 1834.3/4856.879 & 1765.9/4783.888 & 1866.1/4839.318 \\
& 0.8 & 1765.5/4727.801 & 1741.4/4802.599 & 1833.3/4852.144 & 1767.2/4797.778 & 1855.3/4859.506 \\
& 0.9 & 1889.9/4757.937 & 1757.0/4821.188 & 1825.3/4852.615 & 1769.7/4823.612 & 1843.0/4842.762 \\
& 1.0 & 2127.1/4765.809 & 2173.0/4820.209 & 2549.8/4596.665 & 2167.0/4825.327 & 2548.5/4594.403 \\
\bottomrule
\end{tabular}
\caption{Comparison of the matching approaches for the 300 vehicle instances, averaged over buffer settings $\mBuffer=1,..,6$.}
\label{tab:matching-cmp-fleet-300}
\end{table}

\begin{table}[tbh!]
\small
\begin{tabular}{llrrrrr}
\toprule
$\omega$ & $\rho$ & $\mMatchingGreedy$ & $\mMatchingWSCGreedy$ & $\mMatchingWSCSRR$ & $\mMatchingWSPGreedy$ & $\mMatchingWSPSRR$ \\
\midrule
\multicolumn{2}{l}{$\omega^1$} & 1532.6/6201.295 & 1701.0/6266.906 & 1745.8/6400.404 & 1654.5/6267.748 & 1669.1/6392.717 \\
\multirow{11}{*}{$\omega^4(\rho)$} & 0.0 & 981.8/6068.406 & 1054.0/6138.428 & 1196.0/6270.627 & 1080.8/6143.080 & 1221.1/6276.299 \\
& 0.1 & 980.9/6072.116 & 1050.9/6138.227 & 1192.8/6273.384 & 1077.8/6147.392 & 1216.0/6278.538 \\
& 0.2 & 981.1/6073.267 & 1049.5/6134.004 & 1187.2/6277.110 & 1075.8/6145.362 & 1213.1/6268.704 \\
& 0.3 & 984.0/6074.576 & 1046.2/6143.305 & 1185.3/6274.331 & 1073.8/6146.142 & 1208.2/6286.012 \\
& 0.4 & 988.2/6080.514 & 1043.7/6142.247 & 1177.8/6274.005 & 1068.0/6157.526 & 1204.4/6284.090 \\
& 0.5 & 1002.2/6083.441 & 1040.2/6152.210 & 1176.8/6271.122 & 1067.7/6150.915 & 1199.1/6281.223 \\
& 0.6 & 1016.8/6099.991 & 1037.2/6166.595 & 1172.0/6281.990 & 1061.4/6171.651 & 1189.7/6290.304 \\
& 0.7 & 1047.0/6119.708 & 1040.5/6163.533 & 1169.5/6277.818 & 1062.5/6170.538 & 1192.6/6286.092 \\
& 0.8 & 1105.9/6130.786 & 1044.2/6184.225 & 1162.2/6290.867 & 1067.0/6184.524 & 1180.4/6306.789 \\
& 0.9 & 1234.7/6162.225 & 1061.9/6215.387 & 1152.1/6293.486 & 1072.4/6226.208 & 1165.1/6296.625 \\
& 1.0 & 1489.3/6165.207 & 1521.1/6206.850 & 1983.5/6069.544 & 1515.6/6212.933 & 1982.3/6064.973 \\
\bottomrule
\end{tabular}
\caption{Comparison of the matching approaches for the 400 vehicle instances, averaged over buffer settings $\mBuffer=1,..,6$.}%
\label{tab:matching-cmp-fleet-400}
\end{table}

\begin{table}[tbh!]
\small
\begin{tabular}{llrrrrr}
\toprule
$\omega$ & $\rho$ & $\mMatchingGreedy$ & $\mMatchingWSCGreedy$ & $\mMatchingWSCSRR$ & $\mMatchingWSPGreedy$ & $\mMatchingWSPSRR$ \\
\midrule
\multicolumn{2}{l}{$\omega^1$} & 946.0/7565.772 & 1118.2/7719.513 & 1217.5/7818.608 & 1077.9/7710.161 & 1144.0/7809.088 \\
\multirow{11}{*}{$\omega^4(\rho)$} & 0.0 & 493.1/7347.679 & 504.0/7432.558 & 691.2/7571.052 & 523.8/7455.587 & 713.6/7563.534 \\
& 0.1 & 492.9/7353.706 & 500.1/7433.571 & 688.4/7562.421 & 522.3/7449.192 & 703.9/7590.416 \\
& 0.2 & 493.2/7351.611 & 496.6/7441.297 & 686.2/7562.654 & 519.1/7446.008 & 705.0/7571.226 \\
& 0.3 & 496.2/7357.245 & 494.5/7443.295 & 679.6/7569.114 & 518.4/7452.287 & 696.0/7584.730 \\
& 0.4 & 499.7/7360.454 & 493.1/7436.825 & 670.0/7585.827 & 513.5/7451.000 & 692.3/7595.029 \\
& 0.5 & 510.7/7367.604 & 492.4/7443.115 & 670.1/7579.171 & 511.9/7453.013 & 688.7/7575.356 \\
& 0.6 & 519.8/7391.090 & 489.3/7452.531 & 664.7/7576.530 & 508.1/7465.410 & 677.1/7592.963 \\
& 0.7 & 545.2/7420.917 & 488.6/7459.393 & 659.2/7573.059 & 506.7/7478.322 & 679.3/7569.343 \\
& 0.8 & 595.3/7447.010 & 495.4/7471.887 & 653.2/7561.762 & 514.2/7481.435 & 663.8/7598.104 \\
& 0.9 & 707.5/7508.623 & 514.2/7495.982 & 637.5/7573.198 & 524.4/7505.240 & 648.0/7567.485 \\
& 1.0 & 939.2/7569.863 & 935.4/7579.394 & 1518.0/7429.612 & 930.9/7586.929 & 1515.7/7430.132 \\
\bottomrule
\end{tabular}
\caption{Comparison of the matching approaches for the 500 vehicle instances, averaged over buffer settings $\mBuffer=1,..,6$.}
\label{tab:matching-cmp-fleet-500}
\end{table}

\begin{table}[tbh!]
\small
\begin{tabular}{llrrrrr}
\toprule
$\omega$ & $\rho$& $\mMatchingGreedy$ & $\mMatchingWSCGreedy$ & $\mMatchingWSCSRR$ & $\mMatchingWSPGreedy$ & $\mMatchingWSPSRR$ \\
\midrule
\multicolumn{2}{l}{$\omega^1$} & 474.0/8792.016 & 643.0/9042.051 & 802.4/9057.959 & 609.9/9001.867 & 728.6/9034.953 \\
\multirow{11}{*}{$\omega^4(\rho)$} & 0.0 & 187.1/8358.231 & 171.4/8255.163 & 343.0/8630.545 & 179.8/8308.890 & 359.7/8637.470 \\
& 0.1 & 187.4/8349.654 & 168.3/8252.061 & 341.7/8633.042 & 178.1/8296.004 & 349.5/8648.709 \\
& 0.2 & 187.6/8346.485 & 165.1/8259.377 & 339.5/8620.907 & 175.4/8299.075 & 347.1/8648.782 \\
& 0.3 & 187.9/8365.431 & 163.4/8262.530 & 331.8/8639.139 & 173.4/8311.725 & 341.1/8650.158 \\
& 0.4 & 190.0/8376.003 & 161.4/8258.889 & 327.7/8608.274 & 170.5/8301.769 & 337.3/8658.075 \\
& 0.5 & 196.2/8384.167 & 160.9/8261.204 & 324.6/8625.020 & 169.5/8298.686 & 335.1/8642.884 \\
& 0.6 & 200.7/8423.686 & 158.7/8269.980 & 322.3/8587.326 & 167.0/8308.498 & 325.8/8631.863 \\
& 0.7 & 213.5/8474.981 & 157.9/8273.689 & 313.7/8600.174 & 165.1/8328.017 & 320.5/8646.648 \\
& 0.8 & 246.8/8527.777 & 157.3/8321.315 & 305.5/8588.429 & 167.7/8343.331 & 311.0/8627.828 \\
& 0.9 & 323.0/8681.856 & 166.5/8379.565 & 287.8/8599.603 & 171.7/8399.199 & 290.0/8622.919 \\
& 1.0 & 520.0/8782.141 & 467.8/8800.594 & 1117.9/8783.759 & 464.5/8810.612 & 1115.1/8785.587 \\
\bottomrule
\end{tabular}
\caption{Comparison of the matching approaches for the 600 vehicle instances, averaged over buffer settings $\mBuffer=1,..,6$.}
\label{tab:matching-cmp-fleet-600}
\end{table}

\begin{table}[tbh!]
\small
\begin{tabular}{llrrrrr}
\toprule
$\omega$ & $\rho$ & $\mMatchingGreedy$ & $\mMatchingWSCGreedy$ & $\mMatchingWSCSRR$ & $\mMatchingWSPGreedy$ & $\mMatchingWSPSRR$ \\
\midrule
\multicolumn{2}{l}{$\omega^1$} & 156.2/9658.758 & 317.5/9999.219 & 479.7/10143.806 & 274.4/9973.591 & 406.9/10126.834 \\
\multirow{11}{*}{$\omega^4(\rho)$}  & 0.0 & 58.8/8745.881 & 46.1/8458.867 & 159.2/9227.587 & 48.9/8519.776 & 157.2/9323.046 \\
& 0.1 & 59.3/8741.495 & 45.7/8448.529 & 158.9/9218.709 & 48.8/8505.500 & 150.6/9312.502 \\
& 0.2 & 59.2/8741.943 & 44.9/8443.738 & 156.8/9215.348 & 48.2/8496.357 & 150.7/9283.886 \\
& 0.3 & 59.0/8760.754 & 44.3/8442.400 & 151.6/9210.082 & 47.6/8495.190 & 145.2/9285.294 \\
& 0.4 & 59.5/8761.770 & 44.0/8430.118 & 145.9/9193.266 & 47.0/8484.818 & 143.4/9292.794 \\
& 0.5 & 60.2/8792.172 & 43.0/8430.900 & 144.3/9201.860 & 45.6/8481.239 & 141.0/9261.781 \\
& 0.6 & 60.3/8846.850 & 41.9/8439.314 & 143.6/9167.630 & 43.5/8495.277 & 133.2/9262.294 \\
& 0.7 & 61.5/8947.845 & 40.3/8447.715 & 134.7/9170.654 & 42.4/8501.484 & 128.7/9268.173 \\
& 0.8 & 68.1/9123.182 & 38.9/8487.215 & 123.7/9177.005 & 40.5/8537.091 & 121.3/9227.817 \\
& 0.9 & 91.8/9485.880 & 39.2/8580.213 & 111.1/9138.746 & 40.6/8610.422 & 108.7/9180.898 \\
& 1.0 & 202.2/9879.019 & 151.2/9669.139 & 788.6/10049.986 & 149.0/9683.299 & 785.5/10056.420 \\
\bottomrule
\end{tabular}
\caption{Comparison of the matching approaches for the 700 vehicle instances, averaged over buffer settings $\mBuffer=1,..,6$.}
\label{tab:matching-cmp-fleet-700}
\end{table}

\begin{table}[tbh!]
\small
\begin{tabular}{llrrrrr}
\toprule
$\omega$ & $\rho$ & $\mMatchingGreedy$ & $\mMatchingWSCGreedy$ & $\mMatchingWSCSRR$ & $\mMatchingWSPGreedy$ & $\mMatchingWSPSRR$ \\
\midrule
\multicolumn{2}{l}{$\omega^1$} & 15.7/9912.701 & 74.8/10720.931 & 281.1/10809.396 & 64.4/10634.283 & 211.2/10865.947 \\
\multirow{11}{*}{$\omega^4(\rho)$} & 0.0 & 15.3/8800.945 & 8.5/8449.284 & 114.0/9186.774 & 8.4/8510.451 & 96.7/9352.123 \\
& 0.1 & 15.3/8801.367 & 8.3/8438.259 & 111.6/9184.487 & 8.2/8501.105 & 90.0/9339.523 \\
& 0.2 & 15.2/8799.070 & 8.2/8428.931 & 110.8/9174.575 & 7.9/8490.644 & 89.5/9317.990 \\
& 0.3 & 15.3/8813.072 & 8.2/8424.803 & 109.1/9144.491 & 7.8/8486.144 & 86.1/9305.730 \\
& 0.4 & 15.2/8823.105 & 7.8/8417.973 & 103.9/9133.980 & 8.0/8474.138 & 85.4/9306.812 \\
& 0.5 & 15.1/8853.028 & 7.6/8418.018 & 100.9/9145.333 & 7.4/8467.719 & 83.6/9274.370 \\
& 0.6 & 15.2/8900.045 & 7.4/8418.666 & 101.4/9104.079 & 6.8/8475.566 & 79.0/9259.791 \\
& 0.7 & 15.1/8994.387 & 7.1/8422.948 & 95.0/9087.233 & 7.0/8477.652 & 76.7/9250.806 \\
& 0.8 & 15.5/9190.426 & 6.4/8456.321 & 87.6/9074.320 & 6.5/8502.463 & 71.7/9202.167 \\
& 0.9 & 17.0/9613.796 & 7.3/8542.373 & 71.3/9069.061 & 7.3/8577.582 & 64.6/9135.124 \\
& 1.0 & 35.2/10432.105 & 16.8/9900.367 & 504.9/11323.463 & 17.8/9906.190 & 502.7/11328.454 \\

\bottomrule
\end{tabular}
\caption{Comparison of the matching approaches for the 800 vehicle instances, averaged over buffer settings $\mBuffer=1,..,6$.}
\label{tab:matching-cmp-fleet-800}
\end{table}

\section{Parameter tuning}
\label{appendix:parameter-tuning-details}

In this section, we detail our parameter tuning process for the metaheuristic approach presented in Section~\ref{section:methodology:metaheuristic}. We conducted incremental experiments in a \gls{ofat} fashion to identify promising parameter values for our resolution approach. We start with initial settings identified during development and preliminary experiments. Then, one by one, we vary the parameters in the following, arbitrary order: i) weights of the sorting criteria or the recreate operator, ii) acceptance criterion, iii) number of iterations in \gls{rnr} ($M_S$ and $M_A$), iv) average number of routes per decomposed problem, v) number of labels per balas-simonetti node (thickness value), vi) ruin parameters $\overline{c}$ and $L$, vii) ruin parameters $\alpha$ and $\beta$, viii) blink rate in recreate operator, ix) perturbation relocate-exchange ratio, x) \gls{kdspp} limits, xi) \gls{kdspp} start time policy, and xii) recreate request limit.
Tables~\ref{tab:appendix:tuning:recreate_weights}--\ref{tab:appendix:tuning:recreate-limit} show the results of the tuning process steps.

Table~\ref{tab:appendix:tuning:final} summarizes the final parameter settings after concluding the tuning.

\begin{table}[!ht]
    \begin{center}
    \small
    \begin{threeparttable}
    \caption{Tuning results for the recreate sorting criteria weights.}
    \label{tab:appendix:tuning:recreate_weights}
    \begin{tabular}{cc}
    \toprule
    Weights & Gap \\
    \midrule
    6,2,1,2,2,2 & 0.163\%/0.180\%\\
    6,2,1,2,2,4 & 0.104\%/0.093\%\\
    6,2,1,2,4,2 & 0.000\%/0.000\%\\
    6,2,1,4,2,2 & -0.758\%/-0.340\%\\
    6,2,4,2,2,2 & -0.312\%/-0.136\%\\
    6,4,1,2,2,2 & 0.208\%/0.152\%\\
    8,2,1,2,2,2 & 0.074\%/0.062\%\\
    \bottomrule
    \end{tabular}
    \begin{tablenotes}[flushleft]
    \footnotesize
    \item \textit{Notes:} Sorting criteria (random, far, close, pickup, start, delivery-end, time window length). We report the gap (unassigned requests / total distance) to the baseline setting (6,2,1,2,4,2).
    \end{tablenotes}
    \end{threeparttable}
    \end{center}
\end{table}

\begin{table}[!ht]
    \centering
    \small
    \begin{threeparttable}
    \begin{tabular}{ccc}
    \toprule
    Method & $T^{\textsc{init}}$ & Gap \\
    \midrule
    Lin.R2R & 0.111 & 1.063\%/0.076\% \\
    Lin.R2R & 0.222 & 0.734\%/0.141\% \\
    Lin.R2R & 0.333 & 0.000\%/0.000\% \\
    Lin.R2R & 0.444 & 0.629\%/0.425\% \\
    Exp.M & 100 & 7.681\%/3.631\% \\
    Exp.M & 200 & 7.097\%/3.251\% \\
    Exp.M & 1000 & 6.228\%/2.686\% \\
    \bottomrule
    \end{tabular}
    \begin{tablenotes}[flushleft]
    \footnotesize
    \item \textit{Notes:} We report the gap (unassigned requests / total distance) to the baseline setting (Lin.R2R 0.333).
    \end{tablenotes}
    \end{threeparttable}
    \caption{Tuning results for the acceptance criteria (Lin.R2R: Linear Record-to-Record; Exp.M: Exponential Metropolis).}
\label{tab:appendix:tuning:aspiration}
\end{table}

\begin{table}[!ht]
    \centering
    \small
    \begin{threeparttable}
    \begin{tabular}{llc}
    \toprule
    $M_A$ & $M_S$ & Gap \\
    \midrule
    20$e^4$ & 10000 & 1.183\%/0.520\% \\
    20$e^4$ & 5000 & 0.434\%/0.099\% \\
    10$e^4$ & 10000 & 1.228\%/0.589\% \\
    10$e^4$ &  5000 & 0.000\%/0.000\% \\
    10$e^4$ &  2500 & 0.284\%/0.174\% \\
    5$e^4$ &  5000 & 1.123\%/0.667\% \\
    5$e^4$ &  2500 & -0.090\%/0.078\% \\
    \bottomrule
    \end{tabular}
    \begin{tablenotes}[flushleft]
    \footnotesize
    \item \textit{Notes:} We report the gap (unassigned requests / total distance) to the baseline setting ($M_A=10e4, M_S=5000$).
    \end{tablenotes}
    \end{threeparttable}
    \caption{Tuning results for the number of iterations in the R\&R ($M_A$) and for each sub-problem ($M_S$).}
\label{tab:appendix:tuning:iterations}
\end{table}

\begin{table}[!ht]
    \begin{center}
    \small
    \begin{threeparttable}
    \caption{Tuning results for the average number of nodes per sub-problem ($\mAvgSplitSize$).}
    \begin{tabular}{cc}
    \toprule
    Range & Gap \\
    \midrule
    $[500,500]$ & -2.008\%/-1.444\%\\
    $[500,1000]$ & -1.229\%/-0.850\%\\
    $[1000,1000]$ & 0.000\%/0.000\%\\
    $[1500,1500]$ & 2.907\%/1.687\%\\
    \bottomrule
    \end{tabular}
    \begin{tablenotes}[flushleft]
    \footnotesize
    \item \textit{Notes:} We report the gap (unassigned requests / total distance) to the baseline setting ([1000,1000]).
    \end{tablenotes}
    \end{threeparttable}
    \label{tab:appendix:tuning:avg-split}
    \end{center}
\end{table}

\begin{table}[!ht]
    \begin{center}
    \small
    \begin{threeparttable}
    \caption{Tuning results for the number of labels maintained per node (thickness) in the Balas-Simonetti neighborhood.}
    \begin{tabular}{cc}
    \toprule
    Thickness & Gap \\
    \midrule
    4 & -0.321\%/-0.132\%\\
    8 & 0.000\%/0.000\%\\
    12 & 0.046\%/0.016\%\\
    16 & 0.107\%/0.059\%\\
    \bottomrule
    \end{tabular}
    \begin{tablenotes}[flushleft]
    \footnotesize
    \item \textit{Notes:} We report the gap (unassigned requests / total distance) to the baseline setting (8).
    \end{tablenotes}
    \end{threeparttable}
    \label{tab:appendix:tuning:bs-thickness}
    \end{center}
\end{table}

\begin{table}[!ht]
    \begin{center}
    \small
    \begin{threeparttable}
    \caption{Tuning results for ruin parameters $\overline{c}$ and $L$.}
    \begin{tabular}{cccc}
    \toprule
     & \multicolumn{3}{c}{$L$} \\ \cmidrule{2-4}
    $\overline{c}$ & 5 & 10 & 20 \\
    \midrule
    5 &  -0.061\%/-0.004\% & -- & -- \\
    10 &  0.046\%/-0.081\% & 0.0\%/0.0\% & 0.046\%/-0.153\% \\
    15 & -0.015\%/-0.118\% & -0.215\%/-0.115\% & 0.614\%/0.477\% \\
    20 & -- & -0.169\%/-0.136\% & -- \\
    \bottomrule
    \end{tabular}
    \begin{tablenotes}[flushleft]
    \footnotesize
    \item \textit{Notes:} We report the gap (unassigned requests / total distance) to the baseline setting ($\overline{c}=10$, $L=10$).
    \end{tablenotes}
    \end{threeparttable}
    \label{tab:appendix:tuning:ruin-customers-cardinality}
    \end{center}
\end{table}

\begin{table}[!ht]
    \begin{center}
    \small
    \begin{threeparttable}
    \caption{Tuning results for ruin parameters $\mProbRuinMode$ and $\mRuinBeta$.}
    \begin{tabular}{cccc}
    \toprule
     & \multicolumn{3}{c}{$\mProbRuinMode$} \\ \cmidrule{2-4}
    $\mRuinBeta$ & 0.25 & 0.50 & 0.75 \\
    \midrule
    0.01 &  0.200\%/0.119\% &  0.00\%/0.00\% & -0.046\%/-0.032\% \\
    0.10 &  -- &  0.277\%/0.148\% & -0.384\%/-0.242\% \\
    0.50 & -- & 0.200\%/0.110\% & -0.031\%/0.056\% \\
    \bottomrule
    \end{tabular}
    \begin{tablenotes}[flushleft]
    \footnotesize
    \item \textit{Notes:} We report the gap (unassigned requests / total distance) to the baseline setting ($\mProbRuinMode=0.50$, $\mRuinBeta=0.01$).
    \end{tablenotes}
    \end{threeparttable}
    \label{tab:appendix:tuning:ruin-alpha-beta}
    \end{center}
\end{table}

\begin{table}[!ht]
    \begin{center}
    \small
    \begin{threeparttable}
    \caption{Tuning results for the blink rate in the recreate operator.}
    \begin{tabular}{cc}
    \toprule
    Blink & Gap \\
    \midrule
    0.01 & 0.401\%/0.387\%\\
    0.05 & 0.000\%/0.000\%\\
    0.10 &  0.448\%/0.258\%\\
    \bottomrule
    \end{tabular}
    \begin{tablenotes}[flushleft]
    \footnotesize
    \item \textit{Notes:} We report the gap (unassigned requests / total distance) to the baseline setting (0.05).
    \end{tablenotes}
    \end{threeparttable}
    \label{tab:appendix:tuning:blink-rate}
    \end{center}
\end{table}

\begin{table}[!ht]
    \begin{center}
    \small
    \begin{threeparttable}
    \caption{Tuning results for the number of perturbation moves in the ILS, relative to the number of requests in the instance.}
    \begin{tabular}{cc}
    \toprule
    Factor & Gap \\
    \midrule
    1.00 & 0.401\%/0.218\%\\
    1.66 & 0.000\%/0.000\%\\
    2.00 & 0.448\%/0.296\%\\
    \bottomrule
    \end{tabular}
    \begin{tablenotes}[flushleft]
    \footnotesize
    \item \textit{Notes:} We report the gap (unassigned requests / total distance) to the baseline setting (1.66).
    \end{tablenotes}
    \end{threeparttable}
    \label{tab:appendix:tuning:perturbation-moves}
    \end{center}
\end{table}

\begin{table}[!ht]
    \begin{center}
    \small
    \begin{threeparttable}
    \caption{Tuning results for the ratio of relocate moves performed during permutation in the ILS.}
    \begin{tabular}{cc}
    \toprule
    Ratio & Gap \\
    \midrule
    0.25 & 0.370\%/0.214\%\\
    0.50 & 0.000\%/0.000\%\\
    0.75 & 0.015\%/0.054\%\\
    \bottomrule
    \end{tabular}
    \begin{tablenotes}[flushleft]
    \footnotesize
    \item \textit{Notes:} We report the gap (unassigned requests / total distance) to the baseline setting (0.50).
    \end{tablenotes}
    \end{threeparttable}
    \label{tab:appendix:tuning:perturbation-ratio}
    \end{center}
\end{table}

\begin{table}[!ht]
    \centering
    \small
    \begin{threeparttable}
    \begin{tabular}{ccccc}
    \toprule
     & \multicolumn{4}{c}{Time [s]} \\ \cmidrule{2-5}
    Dist. [m] & 600 & 1200 & 1800 & 2400 \\
    \midrule
    3000 & -- & -- & 0.525\%/0.340\% & -- \\
    4000 & 33.812\%/25.852\% & 0.046\%/-0.153\% & 0.0\%/0.0\% & 0.046\%/0.057\% \\
    5000 & -- & -- & 0.448\%/0.197\% & -- \\
    6000 & -- & -- & 0.355\%/0.215\% & -- \\
    \bottomrule
    \end{tabular}
    \begin{tablenotes}[flushleft]
    \footnotesize
    \item \textit{Notes:} We report the gap (unassigned requests / total distance) to the baseline setting (distance: 4000, time: 1800).
    \end{tablenotes}
    \end{threeparttable}
    \caption{Tuning results for KDSP distance and time limits when generating the auxiliary arc.}
    \label{tab:appendix:tuning:kdsp_limits}
\end{table}

\begin{table}[!ht]
    \centering
    \small
    \begin{threeparttable}
    \begin{tabular}{cccc}
    \toprule
    \multicolumn{3}{c}{k-dSPP (start time policy)} \\ \cmidrule{1-3}
    Earliest & Average & Latest & Na\"{i}ve \\
    \midrule
    -0.386\%/-0.132\% & 0.0\%/0.0\% & 0.525\%/0.310\% & 7.623\%/7.566\% \\
    \bottomrule
    \end{tabular}
    \begin{tablenotes}[flushleft]
    \footnotesize
    \item \textit{Notes:} We report the gap (unassigned requests / total distance) to the baseline setting (average).
    \end{tablenotes}
    \end{threeparttable}
    \caption{Tuning results for the recombination methods.}
    \label{tab:appendix:tuning:kdsp_mode}
\end{table}

\begin{table}[!ht]
    \begin{center}
    \small
    \begin{threeparttable}
    \caption{Tuning results for the number of requests considered during recreate.}
    \label{tab:appendix:tuning:recreate-limit}
    \begin{tabular}{cc}
    \toprule
    Limit & Gap \\
    \midrule
    20 & 1.549\%/0.852\%\\
    40 & 0.000\%/0.000\%\\
    60 & 0.186\%/0.078\%\\
    80 & 0.496\%/0.258\%\\
    \bottomrule
    \end{tabular}
    \begin{tablenotes}[flushleft]
    \footnotesize
    \item \textit{Notes:} We report the gap (unassigned requests / total distance) to the baseline setting (40).
    \end{tablenotes}
    \end{threeparttable}
    \end{center}
\end{table}

\begin{table}[!ht]
    \centering
    \small
    \begin{threeparttable}
\begin{tabular}{lllll}
\toprule
Parameter                            & Value                      & \multicolumn{2}{l}{Parameter}           & Value       \\ \cmidrule(l){1-2} \cmidrule(l){3-5}
Acceptance criterion / $T^{\textsc{init}}$ & Lin.R2R / 0.333 & \multicolumn{2}{l}{Ruin parameters} &  \\
ILS perturbation moves $\mPerturbationSize{A}$              & $1.66 \cdot |\mSetRequests|$ &        & $\overline{c}$                 & 15          \\
ILS relocate / exchange move ratio $\mProbRelocateExchange$   & 0.50                       &        & $L$                            & 10          \\
$M_A$, $M_S$                         & 5000, 2500                 &        & $\mProbRuinMode$                       & 0.75        \\
Avg. number of nodes per sub-problem $\mAvgSplitSize$ & 500                        &        & $\mRuinBeta$                        & 0.10        \\
Balas Simonetti Nbh. thickness       & 4                          & \multicolumn{2}{l}{Recreate parameters} &             \\
k-dSPP distance and time limits      & 4000, 1800                 &        & blink rate                     & 0.05        \\
k-dSPP start time policy             & Earliest                   &        & sorting criteria weights       & 6,2,1,4,2,2\\
 & & & request limit & 40\\
\bottomrule
\end{tabular}
    \end{threeparttable}
    \caption{Summary of the final settings found after parameter tuning.}
    \label{tab:appendix:tuning:final}
\end{table}

\section{Fleet minimization}
\label{appendix:fleet-min-algorithm}

The classic \gls{pdptw} benchmark sets use a hierarchical objective function, prioritizing minimizing the fleet size over cost. We integrate a fleet minimization component to compare our approach to such sets and show its efficacy. 
We use a modified \gls{ages}, proposed by \cite{SartoriBuriol2020}, outlined in Algorithm~\ref{alg:ages}.

\begin{algorithm}[!ht]
\caption{Adaptive guided ejection search} 
\label{alg:ages}
\SetKwData{Route}{$\mRoute{}$}\SetKwData{Solution}{$\mSolution$}\SetKwData{CurrentSolution}{$\mSolution'$}
\SetKwData{Stack}{$E$}\SetKwData{PenaltyDecay}{$\mPenaltyDecay$}\SetKwData{PenaltyArray}{$\mPenaltyCounter$}\SetKwData{RequestToInsert}{$u$}
\SetKwData{Counter}{$i$}\SetKwData{MinE}{$min_{\left|E\right|}$}
\SetKwData{MaxPerturbation}{$\mMaxIterations{F}$}\SetKwData{PerturbationSize}{$\mPerturbationSize{F}$}
\SetKwFunction{SelectRandomRoute}{select\_random\_route}
\SetKwFunction{RemoveRoute}{remove\_route}
\SetKwFunction{InitializeStack}{initialize\_stack}
\SetKwFunction{RemoveRequest}{remove\_request}
\SetKwFunction{InsertRequest}{insert\_request}
\SetKwFunction{EjectAndInsert}{eject\_and\_insert}
\SetKwFunction{Perturb}{perturb}
\KwIn{Feasible solution \Solution; Maximum perturbation \MaxPerturbation; Penalty counter \PenaltyArray}
\While{$\Counter < \MaxPerturbation$}{
$\Route \leftarrow \SelectRandomRoute(\Solution)$\;
$\CurrentSolution \leftarrow \RemoveRoute(\Solution, \Route)$\;
$\Stack \leftarrow \InitializeStack(\Route)$\;
$\MinE \leftarrow \left|E\right|$\;
$\PenaltyArray[\RequestToInsert] \leftarrow (1-\PenaltyDecay)\PenaltyArray[\RequestToInsert], \text{for every request } \RequestToInsert$\;\label{alg:ages:update_memory}
    \While{$\Stack \not= \emptyset \vee \Counter < \MaxPerturbation$ }{
        $\RequestToInsert \leftarrow \RemoveRequest(\Stack)$\;
        \If{there is a feasible insertion of~$\RequestToInsert$~in~$\CurrentSolution$}{
            $\CurrentSolution \leftarrow \InsertRequest(\RequestToInsert, \CurrentSolution)$\;
            \If{$\left|\Stack\right| < \MinE$}{ 
                $\MinE \leftarrow \left|\Stack\right|$\;
                $\Counter \leftarrow 0$\;
            }
        }
        \Else{
            $\PenaltyArray[\RequestToInsert] \leftarrow \PenaltyArray[\RequestToInsert] + 1$\;\label{alg:ages:inc_memory}
            $\CurrentSolution \leftarrow \EjectAndInsert(\RequestToInsert, \CurrentSolution, \Stack)$\;
            $\CurrentSolution \leftarrow \Perturb(\CurrentSolution, \PerturbationSize)$\;\label{alg:ages:perturb}
            $\Counter \leftarrow \Counter + 1$\;
        }
    }
    \lIf{$\Stack = \emptyset$}{$\Solution \leftarrow \CurrentSolution$}
}
\Return $\Solution$
\end{algorithm}

It uses ejection search \cite{NagataKobayashi2010}, where an item in a list of unassigned requests is iteratively tested to be inserted into the current solution by potentially removing (ejecting) other requests if necessary. The ejection search is performed sequentially in the number of ejections considered. First, all ejections of a single request are tested. If no feasible insertion is possible, all ejections of two requests are tested. 
To guide the search, a penalty counter, $\mPenaltyCounter$, of failed insertion attempts for each request is maintained during the search and used to prioritize ejections of potentially easy-to-reinsert requests (the fewest number of failed insertion attempts). 
The counter is maintained throughout multiple procedure calls and reinitialized in line~\ref{alg:ages:update_memory} by applying a decay factor $\mPenaltyDecay$ to each entry to nudge the search toward promising regions. 
If an insertion attempt fails, the penalty counter increases in line~\ref{alg:ages:inc_memory}. 
When encountering ties, one possible ejection and insertion opportunity with minimal penalty is randomly selected using reservoir sampling. The sum of the failed insertion attempts is used to guide the search for two ejections. A perturbation procedure at line~\ref{alg:ages:perturb}, right after attempting a potential ejection, aims to diversify the search for the next iteration by performing $\mPerturbationSize{F}$ permutation moves. Herein, we either apply a relocate move with probability $\alpha^{\textsc{ages}}$ or an exchange move otherwise. The search stops when a predefined number of perturbations $\mMaxIterations{F}$ have been attempted. The corresponding tracker is reset if the stack size (the number of unassigned requests) decreases to a new minimum for the given fleet size. 


Our main deviation from existing implementations is introducing the penalty decay factor $\mPenaltyDecay$. In our preliminary experiments, we observed tendencies of the guided search to reach similar prior penalty values before identifying more promising ejections. We avoid some of these time-consuming steps by not completely dismissing prior knowledge. However, we identified a larger factor necessary to prevent inverse effects on the search, where previously hard-to-insert requests need to be considered early to find improvements.

We integrated the described fleet minimization component into our algorithm by running \gls{ages} at the beginning of each iteration of our \gls{ils} approach, defined in Algorithm \ref{alg:decomposition}, right before line \ref{alg:decomposition:decompose}. The penalty counter $\mPenaltyCounter$ is initialized once with a value of 1 for every request and passed to the procedure each iteration. Table~\ref{tab:parameters_ages} lists the parameters used in our experiments.

\begin{table}[!ht]
\begin{tabular}{cccccccc}
\toprule
  \multicolumn{4}{c}{DB-ILS} &
  \multicolumn{4}{c}{AGES} \\ 
  \cmidrule(l){1-4} \cmidrule(l){5-8}
  $\mMaxIterations{A}$ &
  $\mMaxIterationsSplit$ &
  $\mAvgSplitSize$ &
  $\mPerturbationSize{A}$ &
  $\mPenaltyDecay$ &
  $\mMaxIterations{F}$ &
  $\mPerturbationSize{F}$ &
  $\alpha^{\textsc{ages}}$ \\
\midrule
  100000 &
  10000 &
  500 &
  $\left|\mSetRequests\right| \cdot 1.66$ &
  0.8 &
  1000000 &
  10 &
  0.58 \\
\bottomrule
\end{tabular}

\caption{Final parameters settings used in the comparative experiments for the PDPTW benchmark.}
\label{tab:parameters_ages}
\end{table}

\section{Solution structure}\label{appendix:solution-structure}

\begin{figure}
    \centering
    \input{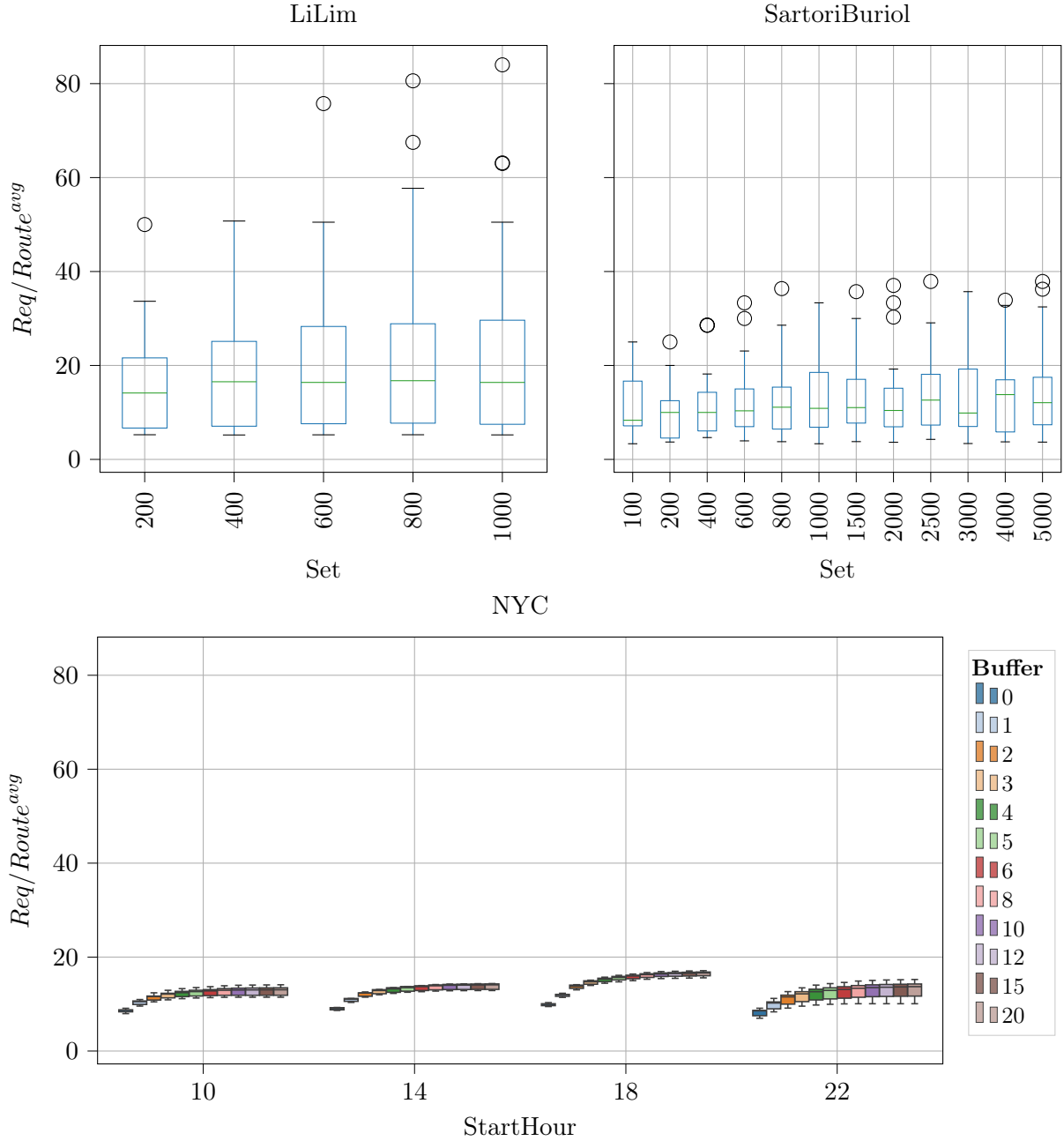}
    \caption{Requests per route for each set in the \cite{LiLim2003}, \cite{SartoriBuriol2020} benchmarks, and the NYC instances.}
    \label{fig:requests_per_route_avg}
\end{figure}

\begin{figure}
    \centering
    \input{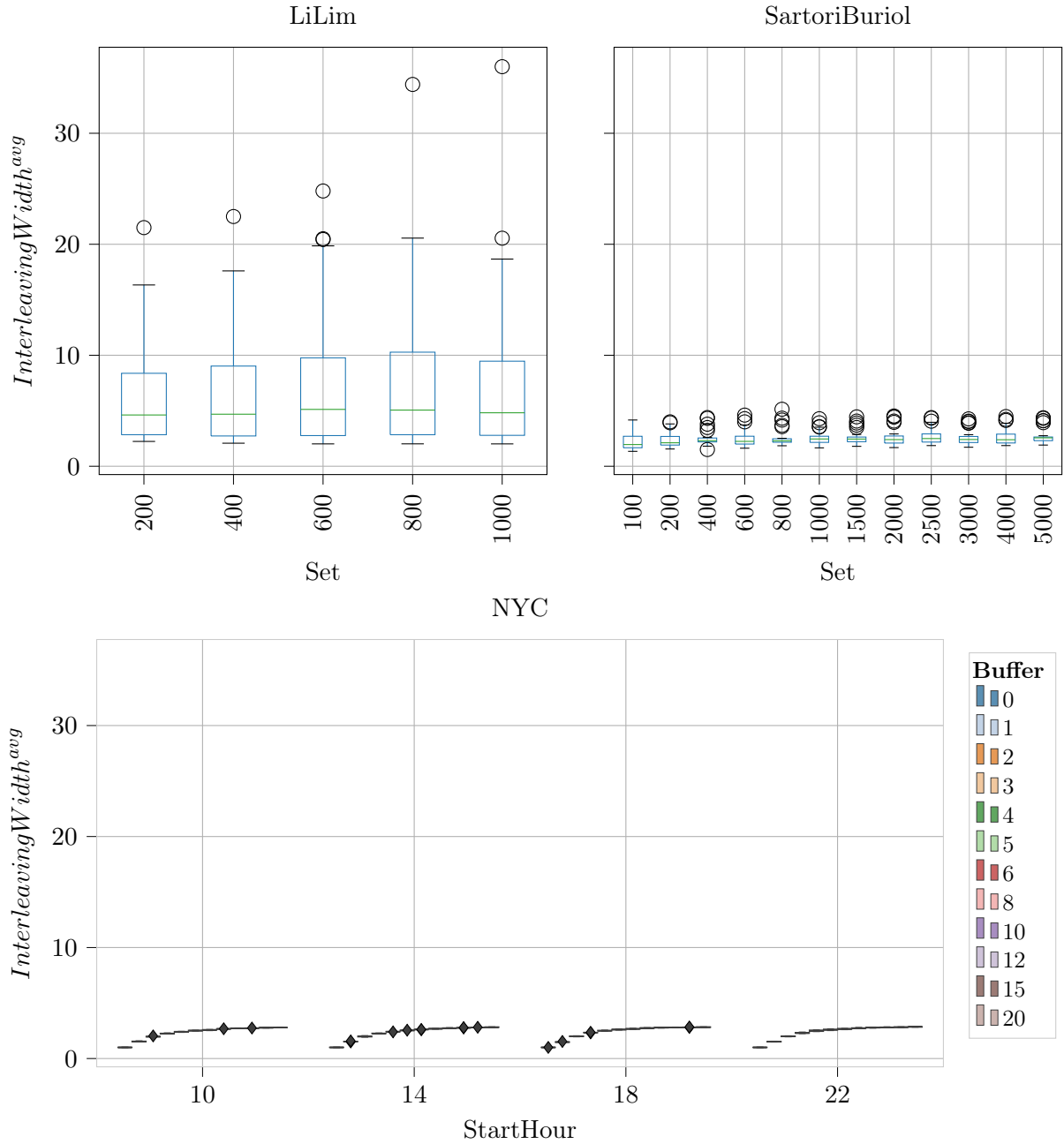}
    \caption{Maximum number of interleaving requests per block for each set in the \cite{LiLim2003}, \cite{SartoriBuriol2020} benchmarks, and the NYC instances.}
    \label{fig:interleaving_width_avg}
\end{figure}

\begin{figure}
    \centering
    \input{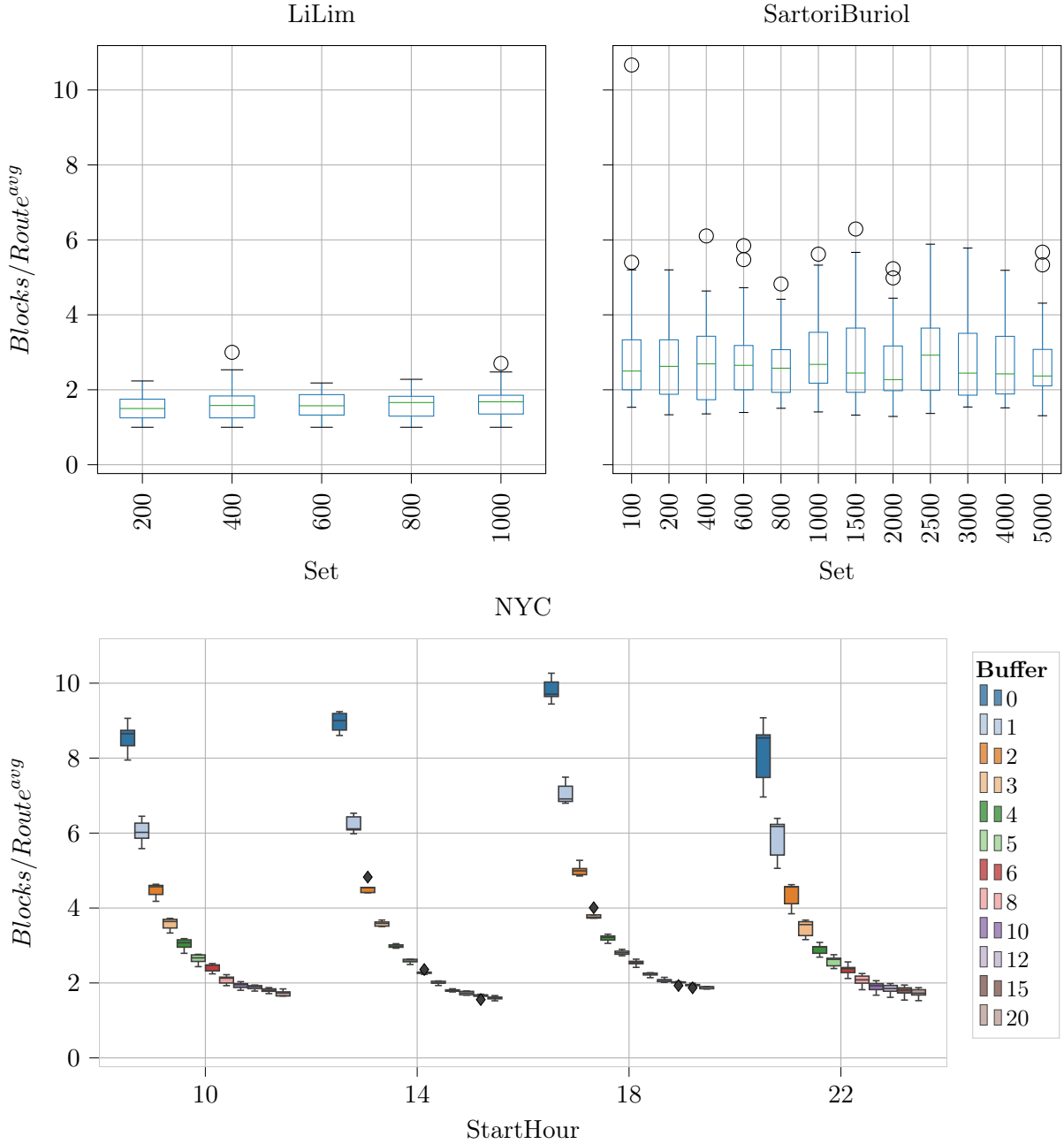}
    \caption{Blocks per route for each set in the \cite{LiLim2003}, \cite{SartoriBuriol2020} benchmarks, and the NYC instances.}
    \label{fig:blocks_per_route_avg}
\end{figure}

With the new set of instances generated based on the NYC Taxi and Limousine Commission \citep{NYC_TLC} data (abbreviated as NYC in this section) and their intended large-scale nature, we set out to identify how good solutions tend to be structured to compare ourselves with the best-fitting benchmark. To analyze the structure of efficient solutions for instances from the literature, we use the best-known solution details from \url{https://www.sintef.no/projectweb/top/pdptw/} and \url{https://github.com/cssartori/pdptw-instances}. In Figures \ref{fig:requests_per_route_avg}--\ref{fig:blocks_per_route_avg}, we compare the average number of requests per route (Figure~\ref{fig:requests_per_route_avg}), the average number of interleaving requests, i.e., how many requests are served concurrently (Figure~\ref{fig:interleaving_width_avg}), and the average number of blocks present per routes (Figure~\ref{fig:blocks_per_route_avg}).

\paragraph{Requests per route.}
We observe in Figure \ref{fig:requests_per_route_avg} that NYC solutions have a low count of requests per route, similar to the \cite{SartoriBuriol2020} benchmark. The solutions for the \cite{LiLim2003} instances tend to have a slightly higher average but spread to up to 8 requests per route. 

\newcommand{\appendixtikzpictureheight}{8cm} 

\paragraph{Interleaving width.}
Figure \ref{fig:interleaving_width_avg} shows a stark difference in solution structure between \cite{LiLim2003}, with many concurrently served requests, and the others, which tend to be similarly low. One reason for this difference lies in the capacity restrictions present. The NYC instances allow up to 3 requests.
\cite{SartoriBuriol2020} instances are generated in two ways: For the New York City instances, they consider passenger transportation with a request demand of up to 6 persons, which is also used as the maximum capacity of vehicles. For the other cities, they consider maximum capacities of between 100 to 300 units of goods, and between 10 to 0.6 times the maximum capacity as demand per request. This leads to an average of 4 to 5.45 requests per vehicle, also reflected in the figure. 
\cite{LiLim2003} instances are based on the classical Solomon instances for the VRPTW and do result in a wider spread with up to more than 30 concurrent requests served.

\paragraph{Blocks per route.}
Figure \ref{fig:blocks_per_route_avg} shows that NYC solutions do consist of more blocks per route (sequences of visits where the vehicle is occupied) with lower buffer settings and gradually move to only two blocks with higher settings. 
Routes in \cite{SartoriBuriol2020} solutions tend to be composed of more blocks (up to more than 6), similar to lower buffer solutions. \cite{LiLim2003} solution routes do average below two blocks, which implies that some vehicles may be occupied from the first visit to a pickup, to their last visit to a delivery.

\paragraph{Summary.}
We observe that solutions of our NYC instances share more similarities with the \cite{SartoriBuriol2020} benchmark set. This informed our decision to focus our design, tuning, and comparison efforts on these instances in our work.



\end{appendices}
\end{document}